\documentclass[11pt]{amsart}

\usepackage{amssymb,latexsym,euscript,amscd,amsmath}

\input xypic

\newcommand{\degr}{\mathrm{deg}}

\theoremstyle{plain}

\newtheorem{thm}{Theorem}
\newtheorem{lem}{Lemma}
\newtheorem{prop}{Proposition}
\newtheorem{cor}{Corollary}

\newtheorem{conjecture}{Conjecture}

\theoremstyle{remark}
\newtheorem{remark}{Remark}
\theoremstyle{definition}
\newtheorem{defn}{Definition}
\title[A conjectural generalization of $n$! result]
{A conjectural generalization of $n$! result
to arbitrary groups}

\author{Shrawan Kumar \and Jesper Funch Thomsen}

\address
{Department of Mathematics, University of North Carolina, Chapel Hill,
NC 27599-3250, USA and Institut for Matematisk Fag, Aarhus Universitet,
Ny Munkegade, DK-8000 \AA rhus C, Denmark.}
\email{kumar$@$math.unc.edu, funch$@$imf.au.dk}

\dedicatory{Dedicated to Prof. M.S. Raghunathan on his sixtieth
birthday}

\begin{document}

\begin{abstract}
The aim of this paper is to formulate a conjecture 
for an arbitrary simple Lie algebra $ \mathfrak g$ in terms of the geometry of 
principal nilpotent pairs.  
 When   $\mathfrak g$ is specialized 
to 
$sl_n$, this conjecture readily implies  
 the $n!$ result and it is very likely that, in fact,  it
  is equivalent to the $n!$ result in this case.  In addition, this 
conjecture can be 
thought of as generalizing an old result of Kostant.
 In another direction, we 
show that to prove the validity of the   $n!$ result for an arbitrary $n$ and an 
arbitrary partition of $n$, it suffices to show its validity only for the 
staircase partitions.
\end{abstract}

\maketitle

\section{Introduction}

To any partition $\sigma : \sigma_0 \geq \sigma_1 \geq \cdots
\geq \sigma_m > 0$ of a positive integer $n$, one associates
 an element $\Delta_\sigma$
of the complex polynomial ring $R_n := \mathbb{C}[X_1,\dots,
X_n,Y_1,\dots,Y_n]$ in $2n$-variables
by defining it (up to a sign) to be the determinant of the matrix
$(X_s^{i_t} Y_s^{j_t})_{1 \leq s,t \leq n}$ for some enumeration 
$\{ (i_1,j_1), \dots, (i_n,j_n) \} $ of the elements in
$D_\sigma := \{ (i,j) \in \mathbb{Z}_+ \times \mathbb{Z}_+ :
i < \sigma_j \}$.  (We could think of  $\Delta_\sigma$ as a 
generalization of  the Vandermonde determinant.) Let  $A_{\sigma}$ be the 
quotient  of the 
ring of constant coefficient differential operators
${\rm Diff}_{2n}$ on $R_n$ by the ideal  $I_{\sigma}$ consisting of the
elements mapping
$\Delta_\sigma$ to zero.

A. Garsia and M. Haiman
\cite{Garsia-Haiman} conjectured that $A_{\sigma}$ should
have complex dimension $n!$. This conjecture later became
known as the {\it $n!$-conjecture}. The interest in the
$n!$-conjecture
came from its connection with the positivity conjecture by
Macdonald, stating that the Kostka-Macdonald coefficients
$K_{\lambda, \mu}(q,t)$ should be polynomials in $q$ and $t$
with nonnegative integer coefficients. Relating the coefficients
of $K_{\lambda, \mu}(q,t)$ to the dimension of certain 
bigraded parts of $A_{\sigma}$, Garsia and Haiman found that
the $n!$-conjecture would imply the above positivity conjecture.

In spite of its simple appearance, the $n!$-conjecture turned
out to be much harder to prove than first anticipated, and
only recently M. Haiman  \cite{haiman} was able to give the
first, and until now only, proof of the conjecture using geometry of the
Hilbert scheme of $n$-points in ${\mathbb C}^2$ (following some suggestions of C. Procesi).
Unfortunately, his proof is quite complicated which makes it
less accessible.

The aim of this paper is to relate the
$n!$-conjecture with the geometry of principal nilpotent pairs ${\bf e}$. More specifically,  let $\mathfrak g$ be a simple 
Lie algebra (corresponding to a
simple adjoint group $G$) with a Cartan subalgebra $\mathfrak h$ and let 
${\bf e} \in 
\mathfrak g^d$ be a 
principal nilpotent pair in $\mathfrak g^d := \mathfrak g \oplus \mathfrak g$  (a
notion defined recently by Ginzburg,
cf. Definition \ref{principal nilpotent pair}). We associate
to ${\bf e}$ a subgroup $\mathcal G^{\bf e}$ of finite index 
of the full automorphism group  $\mathcal G$ of $\mathfrak g$ and consider 
the scheme theoretic intersection of the $\mathcal G^{\bf e}$-orbit 
closure $\bar{\mathcal O}_\mathbf{e}$ of ${\bf e}$ in 
$\mathfrak g^d$ with $\mathfrak h^d$, where  $\mathcal G$ acts on 
$\mathfrak g^d$ diagonally via its adjoint action. Then one of the main
conjectures of this paper asserts that the scheme 
$\bar{\mathcal O}_\mathbf{e}\cap \mathfrak h^d $ is Gorenstein. In addition, 
we conjecture that the affine coordinate ring ${\mathfrak D}_\mathbf{e}:=\mathbb C[\bar{\mathcal O}_\mathbf{e}\cap \mathfrak h^d ]$ supports the regular 
representation of $W$; in particular, it is of dimension $\vert W\vert$ 
(cf. Conjecture \ref{conjecture2}). It is shown (cf. Proposition \ref{our conjecture 
implies n!}) that the above conjecture for $\mathfrak g = sl_n$
readily implies the $n!$ result by using a result of de Concini and Procesi
on the cohomology of Springer fibres and a simple characterization of the algebra  $A_{\sigma}$ in terms of the cohomology of Springer fibres (cf. Theorems 
\ref{characterize}, \ref{haiman algebra}). In fact,  by using the  
deformations of unstable orbits to 
semistable orbits (cf. Section \ref{section4}), we show in  Proposition 
\ref{our conjecture implies n!} that the
(conjectural) Gorenstein 
property of the ring ${\mathfrak D}_\mathbf{e}$ in the case of $\mathfrak g = sl_n$ already implies  
the $n!$ result provided we use a certain variant of 
 a conjecture by Ginzburg 
(cf. Conjecture \ref{conjecture1}). Thus our conjecture provides a 
generalization of  the $n!$ result to an arbitrary simple $\mathfrak g$
($n!$ result being viewed as a result for the case of  the Lie algebra $sl_n$).

In another direction, using the geometry of Hilbert schemes and Borel's fixed 
point theorem,  we show that
the validity of $n!$ result for an arbitrary $n$ and an arbitrary partition of
$n$ follows from its validity for just the staircase partitions 
$\mathfrak s_m : m \geq m-1 \geq \cdots \geq 1 >0$ for any positive integer $m$
(cf. Theorem \ref{staircase}). 

 Let ${\mathbb D}_n^{\mathbb C} (\triangle_n)$ be  the subspace of the polynomial ring $P_n:=
\mathbb C[X_1, \dots, X_n]$ in $n$ variables over $\mathbb C$ obtained by 
applying all the constant coefficient differential operators to  the 
Vandermonde determinant $\Delta_n$. Then it is a well known result that 
${\mathbb D}_n^{\mathbb C} (\triangle_n)$ is of dimension $n!$. Now let $p$ be a 
prime such that $n \leq p^2$ and let $k$ be a field of char. $p$. Then we 
conjecture that the space  ${\mathbb D}_n^k (\triangle_n)$ over $k$ obtained by 
applying all the constant coefficient {\it divided}  differential operators 
to  
 $\Delta_n$ is again  of dimension $n!$ (cf. Conjecture \ref{conjecture3}). This 
conjecture 
immediately implies the validity of $n!$ result for $n=p^2$ and the box 
partition $\mathfrak b_p: p  \geq \cdots \geq p$ of $n$.

For $n = tm+r$ with $t> r \geq 0, t >1$ and $ m \geq 1$, 
consider the {\it box plus one row partition} $\sigma = \sigma(t,m,r) : 
\sigma_0 \geq \sigma_1 \geq \cdots \geq \sigma_m > 0$ of $n$ 
with $\sigma_i = t$ for  $0 \leq i \leq m-1$ and $\sigma_m =r$. 
Finally, in Section \ref{section8}, we realize the algebra $A_\sigma$ for any such partition
as the associated graded
algebra of the cohomology algebra with complex coefficients 
$H^*(SL_n/B, \mathbb C)$
with respect to a certain specific filtration  depending solely on $t$
(cf. Corollary \ref{gr description}), where $B$ is a  Borel subgroup of
$SL_n$ consisting of the upper triangular matrices of determinant one. 
The proof of this corollary relies on a description of the cohomology 
of Springer fibres due to de Concini and Procesi which we have included in the appendix.  In addition, 
 we have included a different  description 
of the cohomology of Springer fibres due to Tanisaki in the 
appendix, which has been used in the proof of Theorem  \ref{haiman algebra}.

We thank V. Ginzburg, M. Haiman,  J.C. Jantzen, N. Mohan Kumar, G. Lusztig, 
and D. Prasad
for some helpful conversations. Part of this work was done while the first 
author was visiting the Newton Institute of Mathematical Sciences, Cambridge during Spring 2001 and the 
Tata Institute of Fundamental Research, Mumbai during Fall 2001. 
The hospitality of these institutions is gratefully acknowledged. 
The first author also acknowledges the support from NSF and one year  
leave from UNC, Chapel Hill. 
The second author was supported by The Danish Research Council, who
also thanks UNC for its hospitality during his visits.

\maketitle
\section{Preliminaries}
\label{section1}
Let $n$ be a positive integer and let $\sigma: \sigma_0\geq \sigma_1\geq
\cdots \geq \sigma_m>0$ be a partition of $n$. The {\it diagram} of 
$\sigma$  is the array of lattice points:
$$D_\sigma := \{(i,j) \in \mathbb Z_+\times \mathbb Z_+: i< \sigma_{j}\},$$
where $\mathbb Z_+ :=\{0,1, \cdots \}$ is the set of nonnegative integers. 
We order the elements of $D_\sigma $ in any manner as $\{1, \cdots, n\}$
and let $ s\in \{1, \cdots, n\}$ correspond to the lattice point 
$(i_s,j_s)$. Define a  generalization of the Vandermonde determinant
as the polynomial in $2n$  variables $X_1, \cdots, X_n, Y_1,
\cdots, Y_n$ with integer coefficients:
$$ \triangle_\sigma := \det \bigl[X_s^{i_t}Y_s^{j_t}\bigr]_{1\leq s,t\leq n}.$$
Clearly a different choice of the ordering of $D_\sigma$ will give rise
to the same $ \triangle_\sigma$ up to a sign. 

Let Diff$_{2n}$ be the algebra ${\mathbb C}[\partial_{X_1}, \cdots, \partial_{X_n}, 
\partial_{Y_1},
\cdots, \partial_{Y_n}]$ of all the constant coefficient differential
operators in $2n$ variables  with coefficients in the field of complex numbers ${\mathbb C}$, where  $\partial_{X_i} := \partial/\partial X_i, 
\partial_{Y_i} := \partial/\partial Y_i$. Let 
Diff$_{2n} (\triangle_\sigma)$ be the space 
$\{D(\triangle_\sigma): D \in  \,\rm{Diff}_{2n}\}. $
The following celebrated conjecture was made by Garsia and Haiman
\cite{Garsia-Haiman}
and proved recently by Haiman \cite{haiman}. 

\medskip

We will refer to the following result 
as $n!$ {\it result for $\sigma$.}

\begin{thm}  ($n!$ result for $\sigma$)  For any positive integer $n$ and partition $\sigma$ of $n$, 
$$ \rm{dim}_\mathbb C \,\rm{Diff}_{2n} (\triangle_\sigma)= n!.$$
\end{thm}

\subsection{Hilbert schemes}

The Hilbert scheme ${\rm Hilb}^n(\mathbb{C}^2)$ of $n$ 
points in $\mathbb{C}^2$ parameterizes the set of ideals 
$I$ in $\mathbb{C}[X,Y]$ whose associated quotient ring 
has complex dimension $n$. By a result of Fogarty 
\cite{Fogarty}, the Hilbert scheme  ${\rm Hilb}^n
(\mathbb{C}^2)$ is a smooth irreducible variety 
of dimension $2n$. Any partition $\sigma$ of $n$ 
defines a point $I_\sigma \in $ Hilb$^n({\mathbb C}^2)$ 
by letting $I_\sigma$ denote the ideal spanned by
the monomials $\{X^iY^j: (i,j) \notin D_\sigma\}.$ 
Thanks to Procesi and Haiman, 
the $n!$ result can be reformulated in geometric 
terms using Hilb$^n({\mathbb C}^2)$ as follows.
 
There is a `tautological' vector bundle ${\mathfrak T}$ 
with the scheme  Hilb$^n({\mathbb C}^2)$ as base, where 
the fibre over any $I\in$  Hilb$^n({\mathbb C}^2)$ is the 
finite dimensional ${\mathbb C}$-vector space (in fact, a 
${\mathbb C}$-algebra) ${\mathbb C}[X,Y]/I$.  
Now consider the  map of sheaves 
$ \mathfrak T^{\otimes n} \otimes \mathfrak T^{\otimes n} \to \wedge^n (\mathfrak T),$
defined by $(a_1\otimes \cdots \otimes a_n) \otimes (b_1\otimes \cdots 
\otimes b_n)\mapsto (a_1b_1)\wedge \cdots \wedge (a_nb_n)$ for any 
$a_i,b_i \in {\mathbb C}[X,Y]/I$ (which is the fibre over $I$).
This map clearly induces the map 
$$\beta: \mathfrak T^{\otimes n} \to  
(\mathfrak T^{\otimes n})^* \otimes \wedge^n (\mathfrak T).$$
Identifying the polynomial ring $R_n:={\mathbb C}[X_1, \cdots, X_n, Y_1,
\cdots, Y_n]$ with the ring ${\mathbb C}[X,Y]^{\otimes n}$, we can view the fibre 
$({\mathbb C}[X,Y]/I)^{\otimes n}$ of the vector bundle 
$\mathfrak T^{\otimes n}$ over $I$  as a quotient of $R_n$.
Let $\bar{K}_\sigma$ be the kernel of $\beta$ over  
$I_\sigma$ and let $K_\sigma$ denote the associated 
preimage in $R_n$. An easy calculation then tells us 
(see \cite{haiman}): 
$$ K_\sigma = \{ f \in R_n : f(\partial_{X_1}, 
\dots, \partial_{X_n}, \partial_{Y_1}, \dots, \partial_{Y_n})(\triangle_\sigma)=0\}.$$
In particular, defining the algebra 
$$A_{\sigma} := R_n/K_{\sigma},$$
we find that $A_\sigma \simeq \rm{Diff}_{2n} (\triangle_\sigma)$ 
as vector spaces. Hence  the $n!$ result is 
equivalent to the statement that the map $\beta$ restricted to the fibre over
$I_\sigma$ has rank $n!$ for all the partitions $\sigma$ 
of $n$. 

It is not difficult to see that the rank of $\beta$ over the 
fibre of any $I \in $ Hilb$_{\rm vg}^n(\mathbb C^2)$ is equal 
to  $n!$, where Hilb$_{\rm vg}^n(\mathbb C^2)\subset $ 
Hilb$^n(\mathbb C^2)$ is the open dense subset consisting 
of those ideals whose zero set consists of $n$ distinct points. 
Using that the rank of $\beta$ along fibres is a lower 
semi-continuous function on ${\rm Hilb}^n(\mathbb{C}^2)$,
we conclude that the rank of $\beta$ along any fibre 
is at most $n!$. Furthermore, any point in ${\rm Hilb}^n
(\mathbb{C}^2)$, by using the natural action of the
torus $(\mathbb{C}^*)^2$ on ${\rm Hilb}^n(\mathbb{C}^2)$,
can be deformed into a point of 
the form $I_\sigma$ for some partition $\sigma$
of $n$.
Applying semi-continuity once more this provides 
us with the following reformulation of the $n!$
result:   

\begin{thm}  
\label{n!}
For any positive integer $n$ and partition $\sigma$ of $n$, 
the rank of $\beta$ over the fibre of $I_\sigma$ is equal to 
$n!.$ Thus, {\rm dim} $A_\sigma = n!$. 

Further, the validity of $n!$ result for all the partitions of $n$ 
is equivalent to the assertion that the rank of $\beta$ over any 
 fibre of  {\rm Hilb}$^n(\mathbb C^2)$  is constant (which is automatically 
equal to $n!$).
\end{thm}

The algebra $A_\sigma$  will play a fundamental role in the paper. Observe that
the algebra $A_\sigma$ is a  bigraded algebra, where the two gradings come respectively from the total degree in the variables $X_1, \cdots, X_n$ and 
the total degree in the variables $Y_1,
\cdots, Y_n$. Moreover, there is a $S_n$-action on $A_\sigma$ 
preserving the bidegrees coming from the (diagonal) action of $S_n$ on
the polynomial ring ${\mathbb C}[X_1, \cdots, X_n, Y_1,
\cdots, Y_n]$, given by: $\theta (X_i) = X_{\theta(i)}$ and
$\theta (Y_i) = Y_{\theta(i)}$, for any $\theta \in S_n$.

We say that a point $I\in$  Hilb$^n({\mathbb C}^2)$ (or alternatively the algebra $\mathbb C[X,Y]/I$; or
alternatively the scheme Spec $\mathbb C[X,Y]/I$) satisfies the {\it 
maximal rank condition} if the rank of $\beta$ over $I$ is equal to 
$n!$ (which is the maximum possible rank over   Hilb$^n({\mathbb C}^2)$).
Thus the set of points   $I\in$  Hilb$^n({\mathbb C}^2)$ satisfying the  
maximal rank condition is open in  Hilb$^n({\mathbb C}^2)$.

\section{A characterization of the algebra $A_\sigma$}
\label{section2}
Fix a positive integer $n$. Let $\sigma:
 \sigma_0\geq \sigma_1\geq
\cdots \geq \sigma_m>0$ be 
a partition of $n$ and let $\sigma^\vee$ be the dual partition. Let
$G=SL_n(\mathbb C)$ and let $X_\sigma$ be the Springer fibre (over complex numbers) corresponding to a fixed  
nilpotent matrix $M_\sigma$ (of size $n\times n$) with Jordan blocks of sizes 
$ \sigma_0,  \sigma_1, 
\cdots,  \sigma_m$. Recall that $X_\sigma$ consists of those Borel 
subalgebras $\mathfrak b$ of $sl_n(\mathbb C)$ such that $M_\sigma \in \mathfrak b$. 
Cohomology of a space $X$ with coefficients in $\mathbb C$ 
will be denoted by $H^*(X,{\mathbb C})$ or simply by $H^*(X)$. In this paper, we
only encounter spaces $X$ such that  $H^*(X)$ is concentrated 
in even degrees. {\it We will consider them as graded algebras under rescaled grading by assigning degree $i$ to the elements of  $H^{2i}(X)$. }

For the definition of Gorenstein rings and their general properties, the 
reader is referred to, e.g., \cite{Bass}.

\begin{thm} 
\label{characterize}
For any partition $\sigma$, there exists an algebra
$T_\sigma$ satisfying the following:

\begin{enumerate}
\item $T_\sigma$ is a $S_n$-equivariant graded quotient of 
$H^*(X_\sigma \times 
X_{\sigma^\vee})$, under the diagonal Springer action of $S_n$ on
$H^*(X_\sigma \times 
X_{\sigma^\vee})$.
\item $T_\sigma^{d_\sigma}$ is one dimensional and, moreover, it is the 
sign representation of $S_n$, where $T_\sigma^m$ denotes the 
graded component of degree $m$ of $T_\sigma$ and $d_\sigma$ is the complex
 dimension of 
the variety $X_\sigma \times 
X_{\sigma^\vee}$. 
\item $T_\sigma$ is a Gorenstein algebra. (Since $T_\sigma$ is a finite 
dimensional graded algebra and its top graded component is one dimensional, 
the condition that $T_\sigma$ is Gorenstein is equivalent to the condition 
that  the pairing
$$ T_\sigma^m \times T_\sigma^{d_\sigma - m} \to T_\sigma^{d_\sigma}$$
induced by the multiplication map is perfect for all $m\geq 0$.) 
\end{enumerate}

Further, such a  $T_\sigma$ (satisfying the above three properties)
is unique up to an isomorphism. More specifically, if
$\varphi_\sigma: H^*(X_\sigma \times 
X_{\sigma^\vee})\to T_\sigma$ and $\varphi'_\sigma: H^*(X_\sigma \times 
X_{\sigma^\vee})\to T'_\sigma$ are two quotients satisfying the above properties,
then their kernels are, in fact, equal. 
\end{thm}

\begin{proof} By the Springer correspondence (see, e.g.,
\cite{Carter}, 13.3 or \cite{Lusztig}, $\S$4.4), the top cohomologies of $X_\sigma$ and 
$X_{\sigma^{\vee}}$ are 
irreducible representations of $S_n$ which differ by
the sign representation $\epsilon$. Hence,
there is a unique copy of $\epsilon$  
in $S_\sigma^{d_\sigma}$, where we abbreviate $H^*(X_\sigma \times 
X_{\sigma^\vee})$ by $S_\sigma$. Let  $(S_\sigma^{d_\sigma})^\epsilon$
be the corresponding isotypical component, which is a one dimensional
space. Now, consider the symmetric bilinear form $\gamma$ obtained as the composition
$$  S_\sigma \times S_\sigma \to S_\sigma^{d_\sigma} \to 
(S_\sigma^{d_\sigma})^\epsilon \simeq \mathbb C\,,$$
where the first map is induced by the multiplication and the second map is 
induced by the $S_n$ equivariant projection. Now define 
$$ T_\sigma := S_\sigma/ {\rm rad}\, \gamma\,,$$
where $ {\rm rad}\, \gamma$ is the radical of $\gamma$. 
It is easy to see that  $ T_\sigma $ satisfies all the three properties listed 
above.

To prove the uniqueness, let $J:= \,{\rm ker}\, \varphi_\sigma$ and
  $J':= \,{\rm ker}\, \varphi'_\sigma$. Then if $J'$ is not a subset of $J$, take 
a nonzero homogeneous element $x\in J+J'/J \subset T_\sigma$. Then, since
$T_\sigma$ is Gorenstein, there exists an element $y\in T_\sigma$ such that
$xy \neq 0$ in $T_\sigma^{d_\sigma}$. We can take a homogeneous 
preimage $\bar{x}$ of $x$ in $J'\subset S_\sigma$. But $J'$ being an ideal,
$\bar{x}\bar{y} \in J'$ for any $\bar{y} \in S_\sigma$. Since
$(J'\cap  S_\sigma^{d_\sigma})^\epsilon = 0$, this leads to a contradiction,
proving that $J'\subset J$. Similarly,  $J\subset J'$, and hence $J=J'$.
\end{proof}

\begin{remark}
\label{manifold}
(a) It is easy to see that  the K{\" u}nneth decomposition of  
$H^*(X_\sigma \times 
X_{\sigma^\vee})$ gives rise to a bigrading on
 $T_\sigma$  and this bigrading respects the $S_n$-action. 

(b) It is well known that $H^*(X_\sigma)$ is concentrated in even degrees
(cf. \cite{DLP}, Theorem 3.9). Moreover, it  can be easily shown that $T^1_\sigma = H^2(X_\sigma \times 
X_{\sigma^\vee})$.

(c) It is natural to ask if there exists a compact topological 
manifold $M_\sigma$  of (real) dimension $2 d_\sigma$ admitting a $S_n$-action on $H^*(M_\sigma) $
and a continuous map $\phi: M_\sigma \to X_\sigma \times 
X_{\sigma^\vee}$ such   that the induced map $\phi^*$ in cohomology is
$S_n$-equivariant and surjective. If so, $H^*(M_\sigma)\simeq T_\sigma $.

We give such a construction for the partition $\sigma: 2 \geq 1$ of 
$n=3$:\footnote{R. MacPherson has mentioned to us that he also constructed such a 
(possibly different) $M_\sigma$ in this case.}   

In this case $\sigma=\sigma^\vee$ and $X_\sigma$ is nothing but two copies of $\mathbb P^1$ 
identified at $\infty$. Thus  $X_\sigma \times 
X_{\sigma^\vee}$ can be identified with the subset  $(\mathbb P^1 \times \infty \times 
\mathbb P^1 \times \infty ) \cup (\infty
\times \mathbb P^1 \times \infty \times \mathbb P^1 )\cup 
(\mathbb P^1 \times \infty \times \infty
\times \mathbb P^1) \cup (\infty \times \mathbb P^1 \times \mathbb P^1 \times \infty )$ of $(\mathbb P^1)^4$. 
Now
take the union of the last two subsets which is join of two copies of $\mathbb P^1
\times \mathbb P^1$ along a point. Finally, take $M_\sigma$ to be the manifold-join of these two
(i.e.,  remove a small disc from these two separately and join them along
the boundary). Then $H^*(M_\sigma)$ is a surjective image of that of 
$H^*(X_\sigma \times X_\sigma^\vee)$. The
Weyl group $S_n$ acts
 on
$H^*(M_\sigma)$ making the above map $S_n$-equivariant (but $S_n$ does not 
act in any natural manner on $M_\sigma$ itself).

(d) It is well known (see, e.g., \cite{Steinberg}, Section 5) 
that $d_\sigma = \sum_{(i,j) \in D_\sigma}\,(i+j),$
where $D_\sigma$ is the diagram of $\sigma$ as in 
Section \ref{section1}. 
\end{remark}

\begin{thm} 
\label{haiman algebra} For any partition $\sigma$ of $n$,  there exists 
a bigraded $S_n$-equivariant isomorphism of algebras: 
 $$T_\sigma \simeq A_\sigma,$$
where $A_\sigma$ is defined in Section \ref{section1}. 
\end{thm}

\vskip1ex
As a preparation to prove the above theorem, we prove the following two lemmas.

Let $\sigma : \sigma_0 \geq
\sigma_1 \geq \cdots \geq \sigma_m >0 $ be a partition of $n$ and let 
 $\sigma^\vee :
\sigma'_0 \geq \cdots \geq \sigma'_{m'}> 0 $ be the dual partition. For any
$m' < i \leq n$, set $\sigma'_i=0,$  and  define the 
integer (for any $1 \leq k \leq n$):
\begin{equation}
\label{dk}
d_k(\sigma) := n - \sum_{s=0}^{n-k-1} \sigma'_{s}~.~
\end{equation}
Let $S$ be a subset of the diagram $D_{\sigma}$ of cardinality $k$, together with a total ordering
on $S$ written as $S = \{ z_1, \dots, z_k \}$. 
Define $i_s, j_s$ by $z_s = (i_s,j_s)$.
Now define the monomial
$$P(S) = \prod_{s=1}^{k} X_s^{i_s} Y_s^{j_s} \in \mathbb C[X_1,
\dots, X_k,Y_1, \dots,Y_k ].$$
For each integer $1 \leq t \leq k$, let $S_t$ denote the
 subset  
$ \{z_s: s > t \}$  of $S$.

\begin{lem}
\label{lemma1}
Let $S$ be as above and let $r$ be a positive integer satisfying
$k - d_k(\sigma) < r \leq k$.
For any sequence  ${\bf s}:1 \leq s_1 < s_2 < \dots
< s_r \leq k$, define
$$P(S)_{\bf s} = (\partial_
{Y_{s_1}}\partial_{ Y_{s_2}}
\cdots \partial_{ Y_{s_r}})P(S).$$
Then there exist integers $1 \leq t' < t \leq k$
such that $P(S)_{\bf s}$ is symmetric in the variables
$(X_t,Y_t)$ and $(X_{t'},Y_{t'})$ (i.e., it is invariant under the 
involution of $ \mathbb C[X_1,
\dots, X_k,Y_1, \dots,Y_k ]$ given by $X_t \mapsto X_{t'}, Y_t \mapsto Y_{t'}).$
\end{lem}
\begin{proof}
Notice first of all that if $n-k \geq \sigma_0$
then $d_k(\sigma) = 0$ and hence there is no $r$ satisfying $k < r \leq k$. 
In the following we
therefore assume that $n-k < \sigma_0$.
By symmetry, we may assume $s_i =i~,i =1, \dots, r$.
Assume that the statement is not true.

If $j_s = 0$ for an integer $s \leq r$ then clearly
$P(S)_{\bf s} = 0$ and the lemma  follows in this case. Hence, assume that 
$j_s >0$ for all  $s \leq r$.
Consider $(i_s,j_s) \in S_r$  and assume
$(i_s,j_s +1) \in S$. Let $s'$ denote the integer
such that $(i_s,j_{s}+1) = (i_{s'}, j_{s'} )$. If
$s' \leq r$ then the statement follows with $t'=s'$
and $t=s$. Hence we may assume that, for each $(i_s,j_s) \in S_r$,
either $(i_s,j_s+1) \notin S$ or
else $(i_s,j_s +1 ) \in S_r$.

Consider the set
$$M = \{ 0 \leq i < \sigma_0 : (i,j) \in S ~\text{for
all $j = 0,1,\dots, \sigma'_{i}-1$} \}.$$
For each $i \in M$ the element $(i,0)$ is then contained
in $S$ and therefore also in $S_r$. Hence,
$(i,1) \in S_r$ if $(i,1) \in S$,
and then $(i,2) \in S_r$ if $(i,2) \in S$
etc. This way we find that
$(i,j) \in S_r $ for $i \in M$ and any $ 0 \leq j < \sigma'_{i}.$
In particular, the number of elements $r$ in $S \setminus
S_r$ is bounded by
$$ r = | S \setminus S_r | \leq
\sum_{s \notin M, s < \sigma_0} (\sigma'_s-1)
\leq \sum_{s = 0}^{\sigma_0 - |M|-1}
(\sigma'_s-1).$$
As the number $n-k$ of elements in $D_{\sigma} \setminus S$
is bounded below by
$$ n-k = |D_{\sigma} \setminus S|  \geq \sigma_0 - |M|,$$
we conclude
$$  r  \leq \sum_{s = 0}^{n-k-1}
(\sigma'_s - 1) = k - d_k(\sigma),$$
which contradicts the assumption.
\end{proof}

\begin{lem}
\label{lem2} Let $\sigma$ be any partition of $n$. 
Let $k \leq n$ and $r$ be positive integers such that 
 $k -d_k(\sigma) < r \leq k$. Then, for any sequence
${\bf s}:1 \leq s_1 < s_2 < \cdots < s_k \leq n$,
$$e_r(\partial_{Y_{s_1}},
\dots,\partial_{Y_{s_k}})
(\Delta_{\sigma}) = 0\,,$$
where $e_r$ denotes the $r$-th elementary symmetric
polynomial.
\end{lem}
\begin{proof}
By symmetry, we may assume $s_i = i$ for all $1\leq i \leq k$.
Expanding $\Delta_{\sigma}$ as a sum of monomials,
we can write (for some index set $A$)
$$\Delta_{\sigma} = \sum_{\alpha \in A} f_\alpha(X_{k+1},
\dots, X_{n}, Y_{k+1}, \dots, Y_n) P(S_\alpha), $$
where $P(S_\alpha)$ are monomials in the variables
$X_1,\dots,X_k,Y_1,\dots,Y_k$ of the form described
earlier  for some  subsets
$S_\alpha \subseteq D_{\sigma}$ of cardinality $k$ together with some 
total ordering on them.
Now, by Lemma \ref{lemma1}, 
$$  e_r(\partial_{Y_{1}},
\dots, \partial_{Y_{k}})
P(S_{\alpha})$$
is a sum of monomials each of which is symmetric
in two distinct variables $(X_t,Y_t)$ and $(X_{t'},
Y_{t'})$ for $1 \leq t < t' \leq k$. In particular,
the polynomial
$$P : = e_r(\partial_{Y_{1}},
\dots, \partial_{Y_{k}})
(\Delta_{\sigma}) $$
can be written as a sum of monomials in the variables
$X_1,Y_1,\dots,X_n,Y_n$ each of which is symmetric
in the  variables $(X_t,Y_t)$ and
$(X_{t'},Y_{t'})$ for some  $1 \leq t < t' \leq k$. The pair 
$t< t'$ could be different for different monomials occurring in $P$. If $P \neq 0$, fix
a monomial $M$ occurring in $P$ and a pair $t< t'$ such that $M$ is symmetric in
 $(X_t,Y_t)$ and
$(X_{t'},Y_{t'})$. 
Let $\pi \in S_n$ be the  permutation $(t,t')$. As $\Delta_{\sigma}$ is
$S_n$-antisymmetric,  $\pi P = -P$. But $\pi (M) = M$, which is a contradiction. Thus $P=0$. 
\end{proof}

Now we are ready to prove Theorem \ref{haiman algebra}.
\vskip1ex

\begin{proof}
By virtue of Lemma \ref{lem2} and 
Theorem \ref{tanisaki} (of the appendix), we get a $S_n$-equivariant morphism
$$\phi_{\sigma} : H^*(X_{\sigma}) \rightarrow A_{\sigma}$$
induced from the natural map $k[Y_1, \dots,Y_n]
\rightarrow A_{\sigma}$. 

From the definition of $\Delta_\sigma$, it is easy to see that 
$$ \tau \Delta_\sigma = \pm \Delta_{\sigma^\vee}\,,$$
where $\tau$ is the algebra automorphism of $R_n$ taking $X_i \mapsto Y_i, Y_i \mapsto X_i.$ Thus, using the
identification
$$ H^*(X_{\sigma^\vee}) \simeq k[X_1,\dots,X_n]/
J(\sigma^\vee), $$
where $J(\sigma^\vee)$ is the ideal defined in Theorem \ref{tanisaki}, 
we get a morphism
$$\phi_{\sigma^\vee} : H^*(X_{\sigma^\vee}) \rightarrow A_{\sigma}$$
induced from the natural map $k[X_1,\dots,X_n]
\rightarrow A_{\sigma}$. Combining the two, we get a
 surjective graded $S_n$-equivariant
morphism
$$\phi : H^*(X_{\sigma}) \otimes H^*(X_{\sigma^\vee})
\rightarrow A_{\sigma},$$
where  $\phi (y \otimes x):= \phi_{\sigma} (y) \cdot
\phi_{\sigma^\vee}(x)$.

By Remark \ref{manifold} (d),  $A_{\sigma}^{d_\sigma} = \mathbb C \Delta_\sigma.$
Thus,  $A_{\sigma}^{d_\sigma}$ transforms via the sign character of $S_n$. We
next show that  $A_{\sigma}$ is Gorenstein: For any $0 \leq m \leq d_\sigma$, take a nonzero $x\in  A_{\sigma}^{m}$ and take a homogeneous lift $\hat{x}\in R_n.$ Viewing $\hat{x}$ as  a differential operator $D(\hat{x})$ (under the identification Diff$_{2n} \simeq R_n$, cf. Section \ref{section1}), we get  $D(\hat{x}) 
\Delta_\sigma \neq 0$. Now take $\hat{y} \in R_n^{d_\sigma - m} $ such that 
$$ D(\hat{y}) ( D(\hat{x}) \Delta_\sigma) =  D(\hat{y}\hat{x})\Delta_\sigma =1.$$
This shows that the pairing $A_\sigma^m \times A_\sigma^{d_\sigma - m } \to 
A_\sigma^{d_\sigma}$ (obtained by the multiplication) is perfect. Thus, 
$A_\sigma$ satisfies the hypotheses of Theorem \ref{characterize}  and 
hence   the isomorphism  $A_{\sigma} \simeq T_{\sigma}$ follows.
\end{proof} 

\begin{remark}
As pointed out by M. Haiman, the map $\phi_{\sigma}$ (and similarly 
$\phi_{\sigma^\vee}$) defined above is injective. Their
existence (and injectivity) can also be obtained from 
the results in \cite{Bergeron-Garsia} and \cite{Garsia-Procesi}
as explained in \cite{Garsia-Haiman}, Section 3.1. However, for the sake of 
completeness, we have included a complete proof. 
  
\end{remark}

\section{Principal nilpotent pairs and deformation}
\label{section4}
Let $\mathfrak{g}$ be a semisimple Lie algebra over 
${\mathbb C}$, 
and let $G=G_{\rm{ad}}$ be the corresponding semisimple linear 
algebraic group of adjoint type. 
Then, as is well known,
 $G$ can be identified with the
identity component of the group $\operatorname{Aut}
(\mathfrak{g})$ of Lie algebra automorphisms of $\mathfrak{g}$. 
We denote the rank of $\mathfrak g$ by rk($\mathfrak g$).

We let $\mathfrak{g}^d$ denote the direct sum $\mathfrak{g} \oplus
\mathfrak{g}$ and think of $G$ as acting diagonally on 
$\mathfrak{g}^d$ by the adjoint action. The following definition 
is due to Ginzburg. 

\begin{defn}(\cite{Ginzburg}, Definition 1.1)
\label{principal nilpotent pair}
A pair $\mathbf{e}=(e_1,e_2) \in \mathfrak g^d$ is called a {\it principal 
nilpotent pair} if the following three conditions hold:
\begin{enumerate}
\item $[e_1,e_2]=0$. 
\item The common centralizer of $e_1$ and $e_2$ in $\mathfrak{g}$ 
is of dimension $\operatorname{rk}(\mathfrak{g})$.
\item For any pair $t_1,t_2 \in {\mathbb C}^*$, there exists $g \in G$ 
such that 
$$\operatorname{Ad}(g) \mathbf{e} = (t_1 e_1,t_2 e_2).$$
\end{enumerate}

(Observe that by (3), $e_1$ and $e_2$ are ad-nilpotent elements.) 
\end{defn}

We  need the following result due to  Ginzburg.

\begin{thm}(\cite{Ginzburg}, {\rm Theorem} 1.2)
\label{associated semisimple pair} 
Let $\mathbf{e}=(e_1,e_2)$ be a principal nilpotent pair.
Then there exists a pair $\mathbf{h}=(h_1 ,h_2) \in \mathfrak g^d$ 
satisfying: 
\begin{enumerate}
\item $h_1$ and $h_2$ are semisimple.
\item $[h_1,h_2]=0$.
\item $[h_i,e_j]= \delta_{ij} e_j$.
\item The common centralizer of $h_1$ and $h_2$ 
is a Cartan subalgebra in $\mathfrak{g}$.
\item All the eigenvalues of ad$_{h_i}: \mathfrak g \to \mathfrak g$ are integral for
$i=1,2$.
\end{enumerate}
\end{thm}

Such a pair $\mathbf{h}$,  called an  
 {\it associated semisimple pair},  is determined by the principal nilpotent pair
$\mathbf{e}$ uniquely up to a $G$-conjugacy.

Let $\mathbf{e}=(e_1,e_2)$ be a principal nilpotent pair. 
As $e_1$ and $e_2$ are ad-nilpotent and commute, the sum 
$e_1+e_2$ is also ad-nilpotent. For any $t \in {\mathbb C}$ we can 
 form the following automorphism of $\mathfrak{g}$:
$$ g(t) = \operatorname{exp}(t \,\operatorname{ad} (e_1+e_2)),$$
or, alternatively,  consider $g(t)$ as an 
element of $G$. The action of $g(t)$ on the associated 
semisimple pair $\mathbf{h}=(h_1,h_2)$ is then 
described by 
$$ g(t) h_i = h_i - t e_i.$$ 
From this one immediately concludes:

\begin{lem}
All the elements of the form $\mathbf{h}(t) = 
(h_1 - t e_1, h_2 - t e_2)$  lie in the 
$G$-orbit $O_{\mathbf{h}} $ of $\mathbf{h}$.
\end{lem}

 The
 {\it induced cone} $\mathfrak{C}O_\mathbf{h}$ in $\mathfrak{g}^d$ is,
 by definition, equal 
to $\cup_{t \in {\mathbb C}} \,\,t \cdot O_\mathbf{h} \subseteq \mathfrak{g}^d$. By the 
above result, we know that $(t h_1 + e_1, t h_2 + e_2)$ is contained
in $\mathfrak{C}O_\mathbf{h}$ when $t \neq 0$. Hence, $(e_1,e_2)$ is contained in 
the closure $\overline{\mathfrak{C}O_\mathbf{h}}$. 

\medskip

\subsection{Deformation}
For a complex vector space $V$, let ${\mathbb C}[V]$  denote the affine
coordinate ring of $V$. This acquires  a natural grading coming from the identification 
${\mathbb C}[V] = S(V^*)$, where $ S(V^*)$ denotes the symmetric algebra of the vector space dual $V^*$ of $V$. If $X$ is any subset
of $V$, we denote by $I(X)$ the ideal of 
functions in ${\mathbb C}[V]$ vanishing identically on $X$.
The vector space spanned by the  leading terms (i.e., the highest degree components) of all the elements in $I(X)$ 
forms an ideal itself. This ideal is denoted by 
$\operatorname{gr}(I(X))$. For any ideal $I \in {\mathbb C}[V]$, we denote  by 
$\mathcal V (I):=\{x\in V: f(x)=0 $ for all $f\in I\}$
the corresponding closed subvariety of $V$. 

In the following, $H$ will denote a reductive group over 
${\mathbb C}$. We also consider a finite dimensional complex 
representation $V$ of $H$. For a $H$-stable vector subspace 
$M\subset {\mathbb C}[V]$, gr($M$) is again $H$-stable, where 
gr($M$) denotes the vector space spanned by the leading terms
of the elements in $M$. Moreover,
$$ \text{gr}(M) \simeq M, \,\,\,\text{as}\, H\text{-modules}.$$

\begin{defn}
An $H$-orbit $Hv$ in $V$ is called {\it semistable} if
$0$ is not contained in its closure. If $Hv$ is not 
semistable  it is said to be {\it unstable}.  

For any subset  $U \subset V$, we define the {\it induced cone}
 (as earlier)
  $\mathfrak{C}U := 
\cup_{t \in {\mathbb C}} \,\,t \cdot U$, and the {\it associated cone}
 $\mathcal{K}(U)$
as the closed (reduced) subvariety $\mathcal V(\operatorname{gr}(I(U)))$ 
of $V$.

Let $Hv$ be a semistable orbit in $V$. 
It is said that an unstable orbit $Hu$ {\it can be  deformed 
into the semistable orbit} $Hv$, if $Hu$ is contained  
in the  associated cone $\mathcal K({Hv})$ of ${Hv}$ and, moreover, $Hu$ has the
same dimension as $Hv$. 
\end{defn}

 With this notation, we recall the following result
from \cite{borho-kraft}, Satz (3.4). 

 \begin{prop}
\label{ideal identity} Let $V$ be a finite dimensional representation of $H$
and let $Hv$ be a semistable orbit. Then 
$$   \mathcal{K}(Hv) = 
(\overline{\mathfrak{C} U} \setminus \mathfrak{C} 
U) \cup \{0\},$$
where $U:=\overline{Hv}$. 

Moreover, 
if $Hu$ is an unstable $H$-orbit such that  there exists a deformation of $Hu$ into the (semistable) orbit 
$Hv$,  then 
$$I(\overline{Hu}) \supset \operatorname{gr} (I(\overline{Hv}))
.$$
\end{prop}

\subsection{Deforming principal nilpotent pairs} 
We now turn to certain deformations inside 
the finite dimensional representation $\mathfrak{g}^d$ 
of $G$.

\begin{lem}
Let $\mathbf{e}=(e_1,e_2)$ be a principal nilpotent pair. 
Then the $G$-orbit $O_{\mathbf{e}} \subset \mathfrak g^d$ is unstable.   
\end{lem}
\begin{proof}
By Definition \ref{principal nilpotent pair}, each element 
of the form $(t e_1,t e_2)$, $t \in {\mathbb C}^*$, is contained in 
$O_{\mathbf{e}}$. Hence $0$ must be contained 
in the closure.
\end{proof}
\begin{lem} 
\label{closed}
Let $\mathbf{h}=(h_1,h_2)$ be a semisimple pair associated to 
a principal nilpotent pair $\mathbf{e}$.  
Then the  $G$-orbit $O_{\mathbf{h}}$ is closed and hence semistable.  
\end{lem}
\begin{proof}
As $h_1$ and $h_2$ commute, one can find a Cartan subalgebra
of $\mathfrak{g}$ containing both $h_1$ and $h_2$. Hence, the common 
centralizer of $h_1$ and $h_2$ in $G$ contains a maximal torus. 
This implies that the orbit is closed.
\end{proof}

\begin{lem}
\label{deformed}
The orbit $O_{\mathbf{e}}$ can be deformed into $O_{\mathbf{h}}$.  
\end{lem}
\begin{proof}
As observed earlier,  the element $\mathbf{e}$ 
lies in $\overline{\mathfrak{C} O_\mathbf{h}}$ and hence, by Proposition
\ref{ideal identity},  also lies
in the  associated cone $\mathcal{K} (O_\mathbf{h})$ 
of $O_\mathbf{h}$ (since all the elements of 
$\mathfrak{C} O_\mathbf{h}$ are semisimple pairs). It remains to 
show that the dimensions of $O_{\mathbf{e}}$ and $O_\mathbf{h}$ 
are the same. For this it is enough to show that the 
dimension $d_\mathbf{e}$ of the common centralizer of $e_1,e_2$ in $\mathfrak g$ 
is  the same  as the dimension   $d_\mathbf{h}$ of the common centralizer 
 of $h_1,h_2$ in $\mathfrak g$. But, by definition, 
$d_\mathbf{e} = d_\mathbf{h} = \operatorname{rk}(\mathfrak{g})$. 
\end{proof}

\subsection{Orbits of the automorphism group of $\mathfrak g$} 
\label{section 4.3}

For the rest of this section we  assume,  for simplicity,  
that $\mathfrak g$ is simple (and, as earlier, 
 $G$ is the associated adjoint group). We 
 fix a maximal torus $T$
of $G$ and  a set of simple roots $\Delta$. The full automorphism 
group $\mathcal G$ of the Lie algebra $\mathfrak g$ is then isomorphic to the 
semidirect product of $G$ with the group $\Gamma$  of the 
diagram automorphisms 
of $\mathfrak{g}$ (for the Dynkin diagram of $\mathfrak g$ with respect to
$\Delta$).

When $e$ is a nilpotent element in $\mathfrak g$, by $\mathcal G^e$ we mean
the subgroup of those elements 
in $\mathcal G$ leaving the $G$-orbit of $e$ invariant. In 
particular, $G \subseteq \mathcal G^e \subseteq \mathcal G$. 
The following lemma gives a complete description of $\mathcal G^e$.

\begin{lem}
\label{G^e}
When $\mathfrak g$ is of type $D_l$ with  even $l>4$, and
$e$ corresponds to a very even partition (with the usual 
classification of nilpotent classes  \cite{Carter}), then $\mathcal G^e$ 
coincides with $G$. When $\mathfrak g$ is of type $D_4$ 
and $e$ corresponds to a very even partition or one of the
partitions $(1^35)$ or $(1^53)$, then $\mathcal G^e$ is a
subgroup of $\mathcal G$ of index 3. In all the remaining 
cases $\mathcal G^e$ coincides with $\mathcal G$. 
\end{lem}
\begin{proof}
Up to $G$-conjugacy we can embed $e$ into a standard 
$sl_2$-triple $(h,e,f)$ with $h$ contained in the 
Lie algebra $\mathfrak h$ of $T$ and such that  
$\alpha(h)$ is positive for all $\alpha \in \Delta$. 
The weighted Dynkin 
diagram of $e$ is then the Dynkin diagram of $\mathfrak 
g$ with respect to $\Delta$, with values $\alpha(h)$, 
$\alpha \in \Delta$, attached to the vertices. 

Fix a diagram automorphism $g$ and consider the induced 
$sl_2$-triple $(g \cdot h, g \cdot e , g \cdot f)$. 
Then $\alpha(g \cdot h)$ is positive for all $\alpha \in \Delta$,
and hence the associated 
weighted Dynkin diagram of $g \cdot e$ is got by applying 
the diagram automorphism $g$ to the weighted Dynkin diagram of 
$e$. As the weighted Dynkin diagram of $e$ determines the 
$G$-conjugacy class of $e$, we see that $g \in \mathcal G^e$ 
if and only if the weighted Dynkin diagram of $e$ is invariant 
under $g$.

We thus conclude that $\mathcal G^e/G$ is the subgroup of $\Gamma$ consisting
of those
diagram automorphisms of the Dynkin diagram of
$\mathfrak g$ which leave the weighted Dynkin 
diagram of $e$ invariant. Now, the result follows by the
inspection of the possible weighted Dynkin diagrams
as found in \cite{Carter}.
\end{proof} 

Spaltenstein introduced an involution, called the {\it Spaltenstein
involution}, on the class of `special' nilpotent orbits in  $\mathfrak g$ 
(cf. \cite{Spaltenstein}). 

\begin{lem}
\label{dual}
Let $e\in \mathfrak g$ be a nilpotent element such that the orbit
$O_e:=G\cdot e$ is special, and 
let $e'$ denote a nilpotent element in the special nilpotent
orbit dual to $G \cdot e$ under the Spaltenstein involution.
Then, for any $\gamma \in  \mathcal G$,
 $\gamma (O_e)$ is special and, moreover,  $\gamma (O_{e'})$
is dual to the orbit  $\gamma (O_e)$. 

In particular, 
$$ \mathcal G^{e}= \mathcal G^{e'}.$$
\end{lem}
\begin{proof}
The first part is known. However we did not find a precise reference,
so we give a brief outline : The Spaltenstein involution is the 
restriction of a map $d$ 
defined on the set of nilpotent orbits (\cite{Spaltenstein},
Chapter III). When $\mathfrak g$ is of classical type, the 
uniqueness statement Theoreme 1.5, Chap.  III in \cite{Spaltenstein} implies
that $d$ commutes with any automorphism of $\mathfrak g$. 
Combining this with the fact that the range of $d$ coincides
with the special nilpotent orbits, the result follows for 
the classical types. 
When $\mathfrak g$ is of type $E_6$ it follows from the diagrams
in \cite{Carter}, p. 441, that there exists at most one 
order reversing involution on the set of special nilpotent 
orbits. In particular, the Spaltenstein involution must
commute with any automorphism 
of $\mathfrak g$,
as it is known to be order reversing,  and the result follows. In the remaining 
cases $\mathcal G$ coincides with $G$ and the statement is trivial.

Alternatively, as pointed out by G. Lusztig, the first part follows from the 
fact that the  Spaltenstein involution corresponds to tensoring the 
associated Springer representation by the sign  representation. 

To prove the `In particular' statement, take $ \gamma \in\mathcal G^{e}$.
Then $\gamma (O_{e'})$ is dual to $\gamma (O_{e})=O_{e}$, and hence 
  $\gamma (O_{e'})= O_{e'}$. This proves the lemma. 
\end{proof}

\vskip1ex

By \cite{Ginzburg}, Lemma 4.13,  for
a principal nilpotent pair ${\bf e}=(e_1,e_2)$ in 
$\mathfrak g$, the $G$-orbits of $e_1$ and $e_2$ 
are Spaltenstein duals of each other. In
particular, by Lemma \ref{dual}, the groups 
$\mathcal G^{e_1}$ and $\mathcal G^{e_2}$ 
coincide. In the following, we  use the 
notation $\mathcal G^{{\bf e}}$ for this 
common group. 

Now we conjecture the following
variant of a conjecture by Ginzburg
\cite{Ginzburg}, Conjecture (2.5). The original conjecture of Ginzburg is 
false,
as seen in Proposition \ref{counterexample}. 

\vskip1ex

\begin{conjecture}
\label{conjecture1}
Let $\mathbf{e}=(e_1,e_2) \in \mathfrak g^d$ be a principal nilpotent pair. 
Then there is a $G$-module isomorphism:
$$ {\mathbb C}[{\bar{\mathcal O}_{\mathbf{e}}}]\simeq {\mathbb C}[G/T],$$
where $T \subset G$ is a maximal torus and $\mathcal O_{\bf e}$ 
denotes the $\mathcal G^{\bf e}$ orbit of ${\bf e}$ in $ \mathfrak g^d$.
\end{conjecture}

\begin{lem}
\label{invariance}
Let ${\bf e} =(e_1,e_2)$ be a principal 
nilpotent pair with an associated semisimple pair 
${\bf h} =(h_1,h_2)$. Then the ideal ${\rm 
gr}(I(G \cdot {\bf h}))$ is 
invariant under the action of $\mathcal G^{\bf e}$.  
\end{lem}
\begin{proof}
When $\mathcal G^{\bf e}$ coincides with $G$ there 
is nothing to prove. By Lemma \ref{G^e},  this takes care of 
any $\mathfrak g$ of type  $B_l$, 
$C_l$, $E_7$, $E_8$, $G_2$, $F_4$ or $A_1$
and also the cases when $e_1$ (and hence also $e_2$)
corresponds to a very even partition in type 
$D_l$,  even $l>4$.

Suppose now that $\mathfrak g$ is of type 
$A_l$ ($l>1$), $E_6$ or $D_l$ ($l$ odd). Then 
$\mathcal G^{\bf e}/G$ coincides with 
 $ \mathcal G/G$ which is the group of order $2$. By Planche I-IX in
\cite{Bourbaki} it follows that the group 
$\mathcal G/G$ is generated by an automorphism
of $\mathfrak g$, which acts on a chosen 
Cartan subalgebra $\mathfrak h$ of $\mathfrak g$ 
via multiplication by $-1$.  Choosing 
$\mathfrak h$ to contain $h_1$ and $h_2$, since
 the ideal  ${\rm 
gr}(I(G \cdot {\bf h}))$ 
is homogeneous, the result follows in these cases as well.

We are left with the situation when $\mathfrak 
g$ is of type $D_l$, with $l$ even. Notice 
that if $g \in \mathcal G^{\bf e}$ leaves 
the $G$-orbit of ${\bf e}$ invariant, then  
by \cite{Ginzburg}, Theorem 1.2(iv), the $G$-orbit
of ${\bf h}$ is also left invariant by $g$.
We may hence concentrate on the situation 
when $g \in 
\mathcal G^{\bf e}$ does not leave the 
$G$-orbit of ${\bf e}$ invariant. Now 
 referring to the classification of 
principal nilpotent pairs as  in 
\cite{El-Pan}, Proposition 5 and Theorem 2,
it follows that any element  $g \in 
\mathcal G^{\bf e}$  leaves the  
$G$-orbit of ${\bf e}$ invariant 
unless ${\bf e}$ correspond to 
a non-integral rectangular skew-graph
(see \cite{El-Pan} for the notation),
which is the same as saying that $e_1$ 
and $e_2$ correspond to very even 
partitions. Thus, when $l>4$,  this situation is already handled 
above.
 This leaves
us with type $D_4$  and  principal 
nilpotent pairs ${\bf e}=(e_1,e_2)$ 
consisting of elements $e_1$ and $e_2$ 
corresponding to very even partitions. 
In this situation there exists an automorphism 
of $\mathfrak g$ carrying $e_1$ and $e_2$ 
to the nilpotent elements corresponding to 
the partitions $(1^3 5)$ and $(1^5 3)$.
As the partitions $(1^3 5)$ and $(1^5 3)$ 
are not very even, this reduces the problem 
to a situation which has already been handled. 
\end{proof}

\begin{thm} 
\label{gr}
Let $\mathbf{h}$ be a semisimple pair associated to 
a principal nilpotent pair $\mathbf{e}$.  Assume the 
validity of the above Conjecture 
\ref{conjecture1} for  $\mathbf{e}$. Then
$$  I({\bar{\mathcal O}_{\mathbf{e}}}) = \operatorname{gr} 
(I({{ O}_{\mathbf{h}}}))
,$$
where $\mathcal O_{{\bf e}}$ denotes the $\mathcal G^{\bf e}$-orbit 
of ${\bf e}$ and  $O_{\bf h}$ denotes the $G$-orbit of 
${\bf h}$ (which is closed by Lemma  \ref{closed}).
\end{thm}

\begin{proof} 
By Proposition \ref{ideal identity}, Lemma \ref{deformed}
and Lemma \ref{invariance},
$$ I({\bar{\mathcal O}_\mathbf{e}}) \supset  \operatorname{gr}
(I({{O}_\mathbf{h})}),$$
which induces a surjective map of $G$-modules
\begin{equation}
\label{surjective}
\mathbb C[\mathfrak g^d]/ ( \operatorname{gr} (I({{O}_\mathbf{h}})))
 \twoheadrightarrow  \mathbb C[{\bar{\mathcal O}_\mathbf{e}}] \simeq
 \mathbb C[G/T],
\end{equation}
where the last isomorphism follows from Conjecture \ref{conjecture1}.
Moreover,  $ \operatorname{gr}
(I({{O}_\mathbf{h}}))\simeq I({{O}_\mathbf{h}}),$
as $G$-modules, so that the left  side of (2) is isomorphic with 
 $\mathbb C[{{O}_\mathbf{h}}]$ (as $G$-modules). But the $G$-isotropy of
$\mathbf{h}$ can easily be seen to be equal to a maximal
torus $T$ and hence $\mathbb C[{{O}_\mathbf{h}}]
\simeq \mathbb C[G/T]$. This forces the surjective map
(\ref{surjective}) to be an isomorphism, thus  the
 theorem follows. 
\end{proof}

\begin{remark}
Using the classification of principal nilpotent 
pairs by Elashvili and Panyushev (see \cite{El-Pan}
and Appendix to \cite{Ginzburg}), one can work 
out the following description of the irreducible components of 
${\bar{\mathcal O}}_{\bf e}$ for a principal nilpotent 
pair ${\bf e}$.
\begin{description}
\item[$A_l$]: ${\bar{\mathcal O}}_{\bf e}$
is the union of the $G$-orbit closures of ${\bf e}$ 
and its transpose ${\bf e}^t :=(e_1^t, e_2^t).$ These 2 components
coincide precisely when ${\bf e}={\bf e}_\sigma$ for a  
box partition $\sigma$ (see the proof of Proposition \ref{our conjecture implies n!} 
 for the notation
 ${\bf e}_\sigma$).
\item[$B_l,C_l,E_l,G_2,F_4, D_l$]  ($l$ even):  ${\bar{\mathcal O}}_{\bf e}$
coincides with the $G$-orbit closure of ${\bf e}$ (cf. Remark on p. 560 
in \cite{Ginzburg}).
\item[$D_l$] ($l$ odd): If ${\bf e}$ is defined 
by a single skew-diagram (see \cite{El-Pan}, Prop. 5(i)), 
then ${\bar{\mathcal O}}_{\bf e}$
consists of two distinct $G$-orbit closures 
of principal nilpotent pairs (those defined 
by the same skew-graph).  Otherwise,
${\bar{\mathcal O}}_{\bf e}$ coincides with the 
$G$-orbit closure of ${\bf e}$.
\end{description} 
\end{remark}

\section{A conjectural  generalization of $n!$ result to arbitrary groups}

In this section $\mathfrak g$ will denote an arbitrary simple 
Lie algebra over  $\mathbb C$ and $\mathfrak h \subset 
\mathfrak g$ a Cartan subalgebra. Let  $G$ be the
associated complex adjoint group with maximal torus $T$ such that  Lie $T=\mathfrak h$. 
As in the last section, $\mathfrak g^d$ denotes the direct sum 
$\mathfrak g \oplus \mathfrak g$ and also  $\mathfrak h^d:=\mathfrak h \oplus \mathfrak h$. 
Let $W:=N/T$ be the Weyl group of $G$; $N$ being the normalizer of $T$ in $G$. 
For a principal nilpotent pair ${\bf e}$, recall the definition of 
$\mathcal G^{\bf e}$ from the last section.

\begin{prop}
\label{intersection} 
 Let $\mathbf{e}=(e_1,e_2)$ be a principal nilpotent pair
in $\mathfrak g$ and let $\mathcal O_\mathbf{e}$ be the diagonal 
$\mathcal G^{\bf e}$-orbit of 
$\mathbf{e}$ in $\mathfrak g^d$. Consider  the affine coordinate ring 
${\mathfrak D}_\mathbf{e}:=\mathbb C[\bar{\mathcal O}_\mathbf{e}\cap \mathfrak h^d]$,
  where 
 $\bar{\mathcal O}_\mathbf{e}\cap \mathfrak h^d$
denotes the scheme theoretic intersection of $\bar{\mathcal O}_\mathbf{e}$ 
with $\mathfrak h^d$ inside $\mathfrak g^d$. Then
${\mathfrak D}_\mathbf{e}$ is a  bigraded finite dimensional (commutative) algebra,
which admits a natural $W$-module structure compatible with the bigrading.  
\end{prop}

\begin{proof}
Consider the action of the two dimensional torus $\mathbb C^*\times 
\mathbb C^*$ on $ \mathfrak g^d$ by $(z,t). (h_1,h_2) = (zh_1,th_2)$,
for $z,t \in \mathbb C^*$ and $ h_1,h_2 \in \mathfrak g$. By the definition 
of principal nilpotent pairs, $\mathbb C^*\times 
\mathbb C^*$ keeps  $\mathcal O_\mathbf{e}$ (and thus $\bar{\mathcal O}_\mathbf{e}$) stable.
Of course, $\mathbb C^*\times 
\mathbb C^*$ also keeps $\mathfrak h^d$ stable. 
This gives the desired bigrading on  ${\mathfrak D}_\mathbf{e}$. 
 Further, under the diagonal action, 
$N$ keeps $ \mathfrak h^d$ stable and, of course, it keeps 
 $\bar{\mathcal O}_\mathbf{e}$ stable. Thus we get an action of $N$ on 
${\mathfrak D}_\mathbf{e}$. But $T$ acts trivially on  $ \mathfrak h^d $; in particular,
it acts trivially on ${\mathfrak D}_\mathbf{e}$ giving rise to an action of 
$W$ (on ${\mathfrak D}_\mathbf{e}$).

To prove that ${\mathfrak D}_\mathbf{e}$ is finite dimensional, it suffices to prove that
set theoretically the scheme $\bar{\mathcal O}_\mathbf{e}\cap \mathfrak h^d$ is a 
singleton. But this follows immediately since the only element of $\mathfrak g^d$ which 
is both 
nilpotent and semisimple is $0$, where 
 we call an element of  $\mathfrak g^d$  
nilpotent (resp. semisimple) if both of its components are 
nilpotent (resp. semisimple). 
\end{proof}

We come to the following one of the main conjectures of the paper.

\begin{conjecture}
\label{conjecture2}
Let the notation and assumptions be as in the above 
Proposition \ref{intersection}. Assume further that 
$\mathbf{e}$ is non-exceptional (\cite{Ginzburg}, 
Definition 4.1).  Then the scheme $\bar{\mathcal O}_\mathbf{e}\cap 
\mathfrak h^d $ is Gorenstein. 
Moreover, \,\,${\mathfrak D}_\mathbf{e}=\mathbb C[\bar{\mathcal O}_\mathbf{e}\cap \mathfrak h^d ]$ supports the regular 
representation of $W$; in particular, it is of dimension $\vert W\vert.$
\end{conjecture}

\begin{remark}  
(a) When  one of the entries $e_1$ or $e_2$ 
of  ${\bf e} = (e_1,e_2)$ is 
zero, then $\mathcal O_{\bf e}$ coincides with the $G$-orbit 
of ${\bf e}$. Hence, in this case, Conjecture 
\ref{conjecture2} 
reduces to a well known result of Kostant. 
\vskip1ex
(b) An old  preprint of Bezrukavnikov-Ginzburg
 contained a
variant of Conjecture \ref{conjecture2} 
 (cf. \cite{Bezrukavnikov-Ginzburg},  Section 8.3). In addition, M. Haiman has informed us 
that he
also conjectured  around 1993 (though unpublished) a variant of Conjecture
\ref{conjecture2} in the special case of
$\mathfrak g = sl_n$.  
\vskip1ex
(c) In an earlier version of our paper,  Conjecture \ref{conjecture2} was formulated 
with $\bar{\mathcal O}_{\bf e}$ replaced by the $G$-orbit closure
 $\bar{ O}_{\bf e}$ of ${\bf e}$.
This turns out to be  false in general (as pointed out by S. Str\o mme), 
since Ginzburg's original conjecture  \cite{Ginzburg},
Conjecture 2.5 is false (cf. Proposition \ref{counterexample}). 
In the same version we had conjectured that 
$\bar{O}_\mathbf{e}\cap \mathfrak h^d $ is even a complete intersection in $\mathfrak h^d$, but it is false in general (as pointed out by Haiman).

\vskip1ex
(d) By 
computer calculations, using SINGULAR \cite{SINGULAR},  
we have checked the above conjecture for $sl_n, n\leq 4$.
\end{remark}

\begin{prop} 
\label{our conjecture implies n!} 
The validity of the above Conjecture \ref{conjecture2} for $\mathfrak g = sl_n$ 
(and every principal nilpotent pair
$\mathbf{e}$) implies the $n!$ result for any partition of $n$.

Alternatively, assume the validity of Conjecture \ref{conjecture1}
and assume, in addition, that ${\mathfrak D}_\mathbf{e}$ is Gorenstein. 
Then the $n!$ result follows for any partition of $n$.

\end{prop}
\begin{proof} Take a partition $\sigma: \sigma_0\geq \sigma_1\geq
\cdots \geq \sigma_m>0$ of $n$. Let $\{v_1, \cdots, v_n\}$ be the
standard basis of $\mathbb C^n$. Define a nilpotent transformation
$e'=e'_\sigma \in sl_n$ which takes $v_1 \mapsto v_2\mapsto v_3 \mapsto
\cdots \mapsto v_{\sigma_0}\mapsto 0,  v_{\sigma_0+1}\mapsto
 v_{\sigma_0+2}\mapsto  \cdots \mapsto  v_{\sigma_0+\sigma_1}\mapsto 0,
 v_{\sigma_0+\sigma_1+1}\mapsto
 v_{\sigma_0+\sigma_1+2}\mapsto  \cdots \mapsto  
v_{\sigma_0+\sigma_1+\sigma_2}\mapsto 0, \cdots .$
Similarly, define
$e''=e''_\sigma \in sl_n$ which takes $v_1 \mapsto v_{\sigma_0+1}\mapsto
 v_{\sigma_0+\sigma_1+1}\mapsto \cdots \mapsto  v_{\sigma_0+ \cdots +
\sigma_{m-1}+1}\mapsto 0,
v_2 \mapsto  v_{\sigma_0+2}\mapsto \cdots \mapsto v_{\sigma_0+ \cdots +
\sigma_{m-1}+2}\mapsto 0 \cdots .$
Then, by \cite{Ginzburg}, $\mathbf{e}_\sigma:=(e', e'')$
is a principal nilpotent pair (cf. also \cite{neubauer}).  

Consider the maps
\begin{equation}
\notag
\begin{split}
H^*(X_{\sigma^\vee}) \otimes  H^*(X_\sigma) & \simeq
\mathbb C[\overline{G\cdot e'}\cap\mathfrak h] \otimes 
\mathbb C[\overline{G\cdot e''} \cap\mathfrak h] \\
& \simeq \mathbb C[\overline{\mathcal G^{e'} \cdot e'} \cap \mathfrak h] 
\otimes \mathbb C[\overline{\mathcal G^{e''} \cdot e''} 
\cap\mathfrak h] \\
& \simeq \mathbb C[(\overline{(\mathcal G^{{\bf e}_\sigma} \times 
\mathcal G^{{\bf e}_\sigma} )\cdot \mathbf{e}_\sigma})\cap 
\mathfrak h^d] \\
& \twoheadrightarrow {\mathfrak D}_{\mathbf{e}_\sigma},
\end{split}
\end{equation}
where the first isomorphism is due to de Concini and Procesi \cite{ConProc} 
and is bigraded and $W\times W$-equivariant 
and the last map
is surjective (bigraded and $W$-equivariant), being induced by the 
closed immersion 
$$\bar{\mathcal O}_{\mathbf{e_\sigma}}\cap \mathfrak h^d \subset (\overline{(
\mathcal G^{{\bf e}_\sigma} \times \mathcal G^{{\bf e}_\sigma}
)\cdot \mathbf{e}_\sigma})\cap \mathfrak h^d.$$

By our Conjecture \ref{conjecture2}, ${\mathfrak D}_{\mathbf{e}_\sigma}$ is Gorenstein; in particular, its top graded component is one dimensional 
(since its zeroth graded component is one dimensional). Moreover, since it 
supports the regular representation of $W$ (by Conjecture \ref{conjecture2}), the top
graded component will have to 
be the sign representation (being the only representation of dimension one 
apart from the trivial representation and the  trivial representation 
occurs in degree $0$). 
 But since ${\mathfrak D}_{\mathbf{e}_\sigma}$ is a
quotient of $ H^*(X_{\sigma^\vee}) \otimes  H^*(X_\sigma)$,
and the latter has a unique copy of the sign representation which occurs
 in the top degree $d_\sigma$ (cf. \cite{Lusztig}, Section 4.4),  we
conclude that ${\mathfrak D}_{\mathbf{e}_\sigma}^{d_\sigma}$ is nonzero 
(and is  the sign representation of $W$). 
Thus, by Theorems \ref{characterize}, \ref{haiman algebra}, we get that 
$$ {\mathfrak D}_{\mathbf{e}_\sigma} \simeq A_\sigma\,.$$
But by Conjecture \ref{conjecture2} again, dim ${\mathfrak D}_{\mathbf{e}_\sigma}= n!$ and thus 
the $n!$ result follows for the partition $\sigma$ (cf. Theorem \ref{n!}).

Now we  prove the `Alternatively' part of the proposition.  
Let $\mathbf{h}_\sigma$ be a semisimple pair associated to 
$\mathbf{e}_\sigma$. Then, by Theorem \ref{gr},
  $  I(\bar{\mathcal O}_{\mathbf{e}_\sigma}) = \operatorname{gr} 
(I({O_{\mathbf{h}_\sigma}})).$ Thus 
$$  I(\bar{\mathcal O}_{\mathbf{e}_\sigma})+ I(\mathfrak h^d) \subset  
\operatorname{gr} 
(I({O_{\mathbf{h}_\sigma}}) + I(\mathfrak h^d)),$$
 and hence  the $S_n$-module 
$\mathbb C[O_{\mathbf{h}_\sigma}\cap \mathfrak h^d]$ is a quotient of the  
$S_n$-module  ${\mathfrak D}_{\mathbf{e}_\sigma}$.
But the scheme $O_{\mathbf{h}_\sigma}\cap \mathfrak h^d$ is reduced 
and is isomorphic with the $S_n$-orbit  $S_n\cdot \mathbf{h}_\sigma$.
It is easy to see that $S_n$ acts freely on $\mathbf{h}_\sigma$,
thus $\mathbb C[O_{\mathbf{h}_\sigma}\cap \mathfrak h^d]$ is the regular 
representation of $S_n$. In particular, the sign representation occurs
in $\mathbb C[O_{\mathbf{h}_\sigma}\cap \mathfrak h^d]$ and hence
also in ${\mathfrak D}_{\mathbf{e}_\sigma}$. But then, by the  argument 
given above, 
${\mathfrak D}_{\mathbf{e}_\sigma}$ 
has a unique copy of the sign
representation which occurs
 in the top degree $d_\sigma$. Moreover, ${\mathfrak D}_{\mathbf{e}_\sigma}$ 
being Gorenstein (by assumption), ${\mathfrak
D}_{\mathbf{e}_\sigma}^{d_\sigma}$ is one dimensional. 
Thus, by Theorems
\ref{characterize} and 
\ref{haiman algebra}, we again get that 
$ {\mathfrak D}_{\mathbf{e}_\sigma} \simeq A_\sigma ,$ and hence $A_\sigma$
has 
dimension at least $n!$. Finally, since dim $A_\sigma \leq n!$ (from the
semicontinuity argument as in Section 2), we get that dim $A_\sigma =
n!$.  
\end{proof}

The following  was pointed out to us independently by D. Panyushev 
and S. Str\o mme for the  principal nilpotent pair $\mathbf{e}_\sigma$
in $sl_3$ associated to the partition $\sigma: 2\geq 1$.

\begin{prop}
\label{counterexample}
Ginzburg's conjecture \cite{Ginzburg}, Conjecture 2.5, is false in general.
\end{prop}

\begin{proof}  
Let ${\bf e}=\mathbf{e}_\sigma$ denote 
the principal nilpotent pair in $sl_n$ 
corresponding to a partition $\sigma$ of $n$ (as in the proof of
Proposition \ref{our conjecture implies n!}) and let 
${\bf h}$ be  an associated 
semisimple pair. By choosing a suitable 
basis, one can  assume that $-{\bf h}$ is an 
associated semisimple pair of the transposed 
principal nilpotent pair ${\bf e}^t$. But if   Ginzburg's conjecture were true
for both of ${\bf e}=\mathbf{e}_\sigma$ and ${\bf e}^t$, then, by a 
variant of Theorem \ref{gr},   we would 
get that 
$I(\overline{G \cdot {\bf e}}) = 
I(\overline{G \cdot {\bf e}^t})$. This would 
imply that $G \cdot {\bf e}=
G \cdot {\bf e}^t$. This is true however
only  for the principal nilpotent pairs  $\mathbf{e}_\sigma$ defined 
by box partitions $\sigma$. Thus, for any partition $\sigma$ of $n$ which 
is not a box partition,  Ginzburg's conjecture is false. 
\end{proof}

\section{Reduction to the staircase partition}
 
 Let $T \subset \operatorname{GL}_2({\mathbb C})$ be the
maximal torus consisting of the diagonal matrices (of nonzero determinant) and let 
 $B$ be the  Borel 
subgroup consisting of the upper triangular matrices (of nonzero determinant). 
Let  $v_1=(1,0)$ and $v_2=(0,1)$ be the 
 standard basis of  ${\mathbb C}^2$.  The algebra of 
regular functions
on ${\mathbb C}^2$ is identified with ${\mathbb C}[X,Y]$, where $X$ and $Y$ are the 
linear maps determined by 
$$ X(v_1) = 1~,~X(v_2)=0~,~Y(v_1) = 0~,~Y(v_2)=1.$$
Consider the standard representation ${\mathbb C}^2$ of $\operatorname{GL}_2
({\mathbb C})$. This induces a natural representation of $\operatorname{GL}_2
({\mathbb C})$ on ${\mathbb C}[X,Y]$. The 
space of  homogeneous polynomials ${\mathbb C}[X,Y]_d$ of degree $d$ in ${\mathbb C}[X,Y]$ 
forms an irreducible $\operatorname{GL}_2
({\mathbb C})$-submodule   
of dimension $d+1$.  The following lemma is well known.

\begin{lem}
\label{B-repr}
Let $V \subseteq {\mathbb C}[X,Y]_d$ be a nonzero $B$-stable subspace. 
Then there exists an integer $i\geq 0$ such that 
$$V = \operatorname{span}_{{\mathbb C}} \{ Y^d, Y^{d-1} X, \cdots, 
Y^{d-i} X^i \}.$$
Conversely, for any $i\geq 0$,  the vector space 
$$\operatorname{span}_{{\mathbb C}} \{ Y^d, Y^{d-1} X, \cdots 
Y^{d-i} X^i \}$$
is $B$-stable.
\end{lem} 

The action of $\operatorname{GL}_2({\mathbb C})$ on ${\mathbb C}[X,Y]$ clearly gives rise to its 
natural action on 
$\operatorname{Hilb}^n({\mathbb C}^2)$ leaving the {\it punctual
Hilbert scheme} $\operatorname{Hilb}^n_0({\mathbb C}^2)\subset 
\operatorname{Hilb}^n({\mathbb C}^2)$ invariant, 
where $\operatorname{Hilb}^n_0({\mathbb C}^2)$
is the projective subvariety consisting of ideals 
containing some  power of the maximal ideal  $(X,Y)$.

Recall from Section \ref{section1} that any partition $\sigma$ of $n$ gives rise to the element
$I_\sigma \in \operatorname{Hilb}^n({\mathbb C}^2)$. 
Furthermore, it is clear from the definition that $I_\sigma$ lies in the 
punctual Hilbert scheme $\operatorname{Hilb}^n_0({\mathbb C}^2)$. Moreover, the 
$T$-stable points of $\operatorname{Hilb}^n({\mathbb C}^2)$ are precisely the points of the form $I_\sigma$, for some partition $\sigma$ of $n$. Since 
$\operatorname{Hilb}^n_0({\mathbb C}^2)$ is a $\operatorname{GL}_2
({\mathbb C})$-stable projective subvariety, the following result
follows readily from the  Borel's
fixed point Theorem.
 
\begin{lem}
\label{B-fix}
The closure $\overline{B I}$ of the $B$-orbit 
of any $I \in \operatorname{Hilb}^n_0({\mathbb C}^2)$ inside $
 \operatorname{Hilb}^n({\mathbb C}^2)$ 
contains a $B$-fixed point (of the form $I_\sigma$). 
\end{lem}

\begin{lem}
\label{B-invariant}
For a partition $\sigma : \sigma_0 \geq \sigma_1 \geq \cdots \geq \sigma_m > 0$
of $n$, the corresponding element $I_\sigma$ is a $B$-invariant point of 
$\operatorname{Hilb}^n ({\mathbb C}^2)$ if and only if 
$\sigma_i > \sigma_{i+1}$ for all $i=0,1,\cdots, m-1$.

Thus, for a partition $\sigma$ as above such that $I_\sigma$ is $B$-invariant, 
we have either $\sigma_0 > m+1 $
or else $m = \sigma_0 -1$ and $\sigma$ is the staircase partition 
$\mathfrak s_{m +1}: m +1 > m  > \cdots >1$ (i.e., $\sigma_i = 
m+1-i$ 
for $i=0,1,\cdots,m$). 
\end{lem}
\begin{proof}
By Lemma \ref{B-repr},  $I_\sigma$ is $B$-invariant
if and only if
$$ X^i Y^j \in I_\sigma~,~i>0~ \Rightarrow X^{i-1} Y^{j+1} \in I_\sigma.$$
By the definition of $I_\sigma$ this is equivalent to the condition:
$$i \geq \sigma_j ~,~i>0 \Rightarrow i-1 \geq \sigma_{j+1}.$$
This is clearly equivalent to $\sigma_i > 
 \sigma_{i+1}$ for all $i=0,1,\cdots, m-1$. This proves the lemma. 
\end{proof}

\vskip1ex

For $f \in {\mathbb C}[X,Y]$ we let $\operatorname{Hilb}^n({\mathbb C}^2)_f$
denote the subset of $\operatorname{Hilb}^n({\mathbb C}^2)$ defined 
by
$$\operatorname{Hilb}^n({\mathbb C}^2)_f = \{ I \in
\operatorname{Hilb}^n({\mathbb C}^2) : f \notin I \}.$$
We prove
\begin{lem}
\label{open}
$\operatorname{Hilb}^n({\mathbb C}^2)_f$ is an open subset of 
$\operatorname{Hilb}^n({\mathbb C}^2)$.
\end{lem}
\begin{proof} Consider the tautological vector bundle 
 ${\mathfrak T}$ over the base Hilb$^n({\mathbb C}^2)$ of rank $n$ as in Section \ref{section1},
where   the fibre over any $I\in$  Hilb$^n({\mathbb C}^2)$ is the
${\mathbb C}$-algebra ${\mathbb C}[X,Y]/I$. Then the multiplication by 
$f$ on each fibre  induces a map of vector bundles  
$m_f : \mathfrak T \rightarrow \mathfrak T$. The subset of 
$\operatorname{Hilb}^n ({\mathbb C}^2)$ 
where $m_f$ has nonzero rank  is then an open subset 
coinciding with 
$\operatorname{Hilb}^n ({\mathbb C}^2)_f$. 
\end{proof}

For a partition $\sigma : \sigma_0 \geq \sigma_1 \geq \cdots \geq 
\sigma_m > 0$ of $n$,  consider the ideal 
$$I_{\sigma}(\Lambda) = (Y^{m+1}(Y-\Lambda), Y^{m+1} X) +
(X^i Y^j : i \geq \sigma_j ~,~j\leq m )$$
inside the polynomial ring ${\mathbb{C}}[X,Y,\Lambda]$.
The quotient  
${\mathbb{C}}[X,Y,\Lambda]/I_\sigma(\Lambda)$  is then a free 
${\mathbb{C}}[\Lambda]$-module of rank $n+1$ with basis 
$\{ X^i Y^j : i < \sigma_j \}\cup \{Y^{m+1}\}$. In 
particular, this defines  a flat family 
$$\psi : \operatorname{Spec}({\mathbb C}[X,Y,\Lambda]/I_\sigma(\Lambda))
\rightarrow  \operatorname{Spec}({\mathbb C}[\Lambda]) \simeq 
\mathbb{A}^1,$$
of closed subschemes of ${\mathbb C}^2$ of length $n+1$. Notice
that the fibre of $\psi$ over $\lambda \neq 0$ is 
Spec($ {\mathbb C}[X,Y]/ I_\sigma \times {\mathbb C}[X,Y]/(Y-\lambda,X) $),
while the fibre over $\lambda = 0$ is Spec(${\mathbb C}[X,Y]/I_{\sigma'}$)
where $\sigma' : \sigma_0 \geq \sigma_1 \geq \cdots \geq \sigma_m\geq 
\sigma_{m+1} > 0$, with $\sigma_{m+1}=1$, is a partition of $n+1$.

We conclude 

\begin{lem}
\label{extra point}
The validity of the $n!$-result (rather the  $n+1!$-result)  for the partition $\sigma'$ implies 
its validity for the partition $\sigma$, where $\sigma'$ is as above.
\end{lem}
\begin{proof} Since the set of points of  Hilb$^{n+1}({\mathbb C}^2)$ where the maximal
rank condition is satisfied is open  (cf. Section \ref{section1}), we get that it is satisfied 
in an open subset of  Hilb$^{n+1}({\mathbb C}^2)$ containing $I_{\sigma'}$. In particular,
from the above family, we get that the maximal
rank condition is satisfied  for $\psi^{-1}(\lambda)$ for some 
$\lambda \neq 0.$  But $\psi^{-1}(\lambda) \simeq \text{Spec} ( {\mathbb C}[X,Y]/ I_\sigma  \times
{\mathbb C})$.  Finally, observe that the maximal rank condition
is satisfied for the  product of rings $A \times {\mathbb C}$ if and 
only if it is satisfied for $A$. This proves the lemma.
\end{proof}

\begin{thm}
\label{staircase}
To prove the validity of the $n!$-result for all $n\geq 1$ and all
the partitions of $n$, it suffices to prove its validity for only 
the  staircase 
partitions $\mathfrak s_m: m\geq m-1\geq \cdots \geq 1$ (for all $m\geq 1$).  
\end{thm}
\begin{proof}
For any $I\in \operatorname{Hilb}^n_0({\mathbb C}^2)$, we first
define $d=d(I)$ as the smallest integer such that  ${\mathbb C}[X,Y]_{d}
\subseteq I$. (Since $I\in \operatorname{Hilb}^n_0({\mathbb C}^2)$, $d$
exists.)  Then clearly $n \leq  d(d+1)/2$ with equality if and
only if $I=I_{\mathfrak s_d}$ for  the staircase  partition $\mathfrak s_d$.
Now define $N(I)=d(d+1)/2-n$. Thus $N(I)$ measures, in some sense,
how far $I$ is from $I_{\mathfrak s_d}$. The assumption of the theorem
means that the maximal rank condition is valid for any $I$ with
$N(I)=0$. We prove the theorem by induction on $N(I)$.

So assume $N(I_{\sigma}) > 0$ and that the theorem has been
proved for $I_{\mu}$ with $N(I_\mu) < N(I_\sigma)$.

Consider the orbit closure
$\overline{B I_\sigma}$ inside $\operatorname{Hilb}^n({\mathbb C}^2)$.
By Lemma \ref{B-fix}, there exists a $B$-fixed point
$I_\mu \in \overline{B I_\sigma}$ for some  partition
$\mu$. Since, for any $d$, ${\mathbb C}[X,Y]_{d}$ is stable under $B$
(in fact, under GL$_2({\mathbb C})$),  by Lemma \ref{open},
${\mathbb C}[X,Y]_{d(I_\sigma)}$ is contained in $I_\mu$.
In particular, $d(I_\mu) \leq d(I_\sigma)$ and thus 
$N(I_{\mu}) \leq N(I_{\sigma})$. Since the set of points
in Hilb$^n(\mathbb C^2)$ satisfying the maximal rank condition
is open, we may now assume $\sigma = \mu$ or, in other
words, we may assume that $I_\sigma$ is $B$-stable.
Then $\sigma: \sigma_0 >  \sigma_1 > \cdots > \sigma_m
>0 $ from which we conclude that $d(I_\sigma) \geq m+1$
with equality if and only if $\sigma$ is a staircase
partition (use  Lemma \ref{B-invariant}). We  therefore
assume $d(I_{\sigma}) > m+1$.

Now, the partition $\sigma'$ of $n+1$, as defined  earlier,
also satisfies that ${\mathbb C}[X,Y]_{d(I_\sigma)}$
is contained in $I_{\sigma'}$. But $N(I_{\sigma'}) <
N(I_\sigma)$. Thus, by induction,  the maximal rank condition
follows for $I_{\sigma'}$. Finally, by Lemma \ref{extra point},
the maximal rank condition  follows for $I_\sigma$.
\end{proof}

\begin{remark}
Observe that GL$_2 ({\mathbb C})$ has an invariant in   $\operatorname{Hilb}^n({\mathbb C}^2)$
if and only if $n$ is of the form $d(d+1)/2$ and in this case $I_{\mathfrak s_d}$,
for the staircase partition $\mathfrak s_d$, 
is the only invariant. 
\end{remark} 

\vskip2ex

For the staircase partition 
$\mathfrak s_m: m\geq m-1\geq \cdots \geq 1$ of
$n:=m(m+1)/2$, we have the following (stronger) variant of Theorem 
\ref{characterize}.
\begin{prop} Let $m$ be any positive integer and let 
$B_m$ be a finite dimensional graded
$S_n$-equivariant quotient of the polynomial ring $R_n:={\mathbb C}[X_1, \cdots, X_n, Y_1,
\cdots, Y_n]$ such that  $B_m^d =0$ for $d> d_{\mathfrak s_m}$ and $B_m^{d_{\mathfrak s_m}}$
is the one dimensional sign representation of $S_n$, where $n:=m(m+1)/2$.
Assume further that 
$B_m$ is a  Gorenstein algebra.
   Then 
$$B_m \simeq A_{\mathfrak s_m},$$
as graded $S_n$-algebras. 

(Observe that $d_{\mathfrak s_m}= m(m-1)(m+1)/3$.)
\end{prop}
\begin{proof} The proof of this proposition follows by a similar argument as that of Theorem \ref{characterize} in view of the following Lemma \ref{sign}.
\end{proof}

\begin{defn}
Let $n$ be a positive integer. Write $n$ in the form
$$n  = 1 +2 +3 +\cdots + m + s~,~0 \leq s < m+1 .$$
The {\it{degree}} of $n$ is defined as
$$\degr(n) = 1 \cdot 0 + 2 \cdot 1 + \cdots + m \cdot
(m-1) + s \cdot m. $$
The integer $s$ is called the {\it remainder} of $n$.
\end{defn}

\begin{lem}
\label{sign}
The lowest degree in $R_n$ in which the sign representation
appears is exactly $\degr(n)$. If the remainder of $n$ is
$0$ then (and only then) there is a unique copy of the
sign representation in degree $\degr(n)$ in $R_n$.
\end{lem}
\begin{proof}
Assume that $L \subseteq R_n^d$ is a copy of the sign
representation. Let $l \in L$ be  a nonzero element.
Write :
$$l = \sum_{\underline{i}, \underline{j}}
\alpha_{\underline{i},\underline{j}}
{X}^{\underline{i}} {Y}^{\underline{j}},$$
where $\underline{i}=(i_1, \dots, i_n) ,\underline{j}=(j_1,
\dots, j_n) \in \mathbb Z_+^n$, 
$X^{\underline{i}}:= X_1^{i_1}\cdots X_n^{i_n}$ and 
$ Y^{\underline{j}}$ has a similar meaning.
Take a nonzero $\alpha_{\underline{i'},\underline{j'}}$.
If there exist $s \neq t$ such that $(i'_s,j'_s) = (i'_t,j'_t)$
then the permutation $(s,t) \in S_n$ acts on $l$ leaving
the coefficient of ${X}^{\underline{i'}}
{Y}^{\underline{j'}}$ invariant. As $L$ is a copy of 
the sign representation this is impossible. Hence,
$(i'_s,j'_s) \neq (i'_t,j'_t)$ whenever
$s \neq t$, i.e., $(i'_1,j'_1), (i'_2,j'_2), \dots,
(i'_n,j'_n)$ are distinct elements in $\mathbb Z_+^2$.

Consider the function $\phi : \mathbb Z_+^2 \rightarrow \mathbb Z_+$
defined by $\phi(x,y) = x+y$. Call $\phi(x,y)$ the norm of $(x,y)$.
Note that there are exactly $s$ elements in $\mathbb Z_+^2$ of norm
$s-1$.
Notice also that the degree of   ${X}^{\underline{i'}}
{Y}^{\underline{j'}}$ is 
$d = \sum_{t} \phi(i'_t,j'_t).$
As the pairs $(i'_t,j'_t)$ are distinct, we get that $d \geq \degr(n)$.

Now we construct a copy of the sign representation inside
$R_n^{\degr(n)}$: Choose $n$ distinct pairs $(i_t,j_t)$
such that the sum of the norms is $\degr(n)$ (start with
pairs of norm 0, then take those of norm 1, $\cdots$). Notice that there is a unique way of
doing this exactly when the remainder of $n$ is $0$. Put
$$L = {\mathbb C} (\sum_{\sigma \in S_n} \epsilon (\sigma)
\sigma ({X}^{\underline{i}}
{Y}^{\underline{j}})).$$
Then $L$ is a copy of the sign representation sitting in the
degree equal to $\degr(n)$.
\end{proof}

\section{Divided differential operators and the box partition}
Let $n$ be a positive integer and let $R^\mathbb Z_n:=\mathbb Z[X_1, \cdots, X_n, Y_1,
\cdots, Y_n]$ be the polynomial ring over $\mathbb Z$. For any $m\geq 0$ and any $1 \leq s \leq n$, the divided differential operators  $\partial_{X_s}^{(m)} :=
\frac{1}{m!} (\frac{\partial}{\partial
X_s})^m$ and  $\partial_{Y_s}^{(m)} :=
\frac{1}{m!} (\frac{\partial}{\partial
Y_s})^m$ keep $R^\mathbb Z_n$ stable. In particular, for any $\mathbb Z$-module
$M$, we obtain the operators (again denoted by)  $\partial_{X_s}^{(m)}$
and  $\partial_{Y_s}^{(m)}$ on $R^M_n:= R^\mathbb Z_n\otimes_\mathbb Z \,M$. 

{\it Now let $p$ be a prime and $k$ a field of char. $p$.}  We will particularly
be interested in the operators  $\partial_{X_s},  \partial_{X_s}^{(p)},$
 $\partial_{Y_s}$ and  $\partial_{Y_s}^{(p)}$ acting on $R_n^k$. Similarly,
we have the operators  $\partial_{X_s},  \partial_{X_s}^{(p)}$ on $P_n^k:=
k[X_1, \cdots, X_n].$

Define a $k$-linear map (in fact, a $k$-algebra homomorphism)
$\phi:R_n^k \to P_n^k$, by 
$X^{\underline{i}} Y^{\underline{j}} \mapsto X^{\underline{i}+ p
\underline{j}}$, for any $ \underline{i}=(i_1, \cdots, i_n)$ and   
 $ \underline{j}=(j_1, \cdots, j_n)\in \mathbb Z_+^n$. 
\begin{lem}  
\label{divided1} For any $1 \leq s\leq n$ and  $Q \in R^k_n$, 

\noindent
(a) $\phi \partial_{ X_s}=
\partial_{ X_s} \phi$\,, \\
(b) $\phi (\partial_{ Y_s} Q) = \partial_{X_s}^{(p)} \phi (Q)$,
if  deg $Q$ in $X_s$ is $<p$. 
\end{lem}
\begin{proof}    (a) Let $\delta_s := (0, \cdots, 1,0, \cdots, 0)\in \mathbb Z_+^n$, where 
$1$ is placed in the $s$-th place.
$$
\begin{array}{lll}
\phi ( \partial_{ X_s} X^{\underline{i}}
Y^{\underline{j}}) &=& \phi (i_s (X^{\underline{i} - \delta_s }Y^{\underline{j}})) \\
&=&i_s X^{\underline{i} - \delta_s + p \underline{j}}.
\end{array}$$
Also,
$$
\begin{array}{lll}
 \partial_{ X_s} \phi
(X^{\underline{i}}Y^{\underline{j}}) &=&  \partial_{ X_s}
X^{\underline{i}+p \underline{j}} \\
&=& i_s X^{\underline{i}+p \underline{j} - \delta_s}. 
\end{array}
$$
This proves (a). 

(b) Let $Q=X^{\underline{i}} Y^{ \underline{j}}$. Then
$$
\begin{array}{lll}
\phi (\partial_{Y_s} Q)&=& \phi (j_s
X^{\underline{i}} Y^{\underline{j} - \delta_s})\\
&=& j_s X^{\underline{i}+p (\underline{j} - \delta_s)}.
\end{array} $$
Also,
$$
\begin{array}{lll}
\partial_{X_s}^{(p)} \phi (Q) &=& \partial_{X_s}^{(p)}
X^{\underline{i}+p \underline{j}} \\
&=& \binom{i_s+pj_s} {p} X^{\underline{i}+p \underline{j} -p\delta_s} 
\\
&=& j_s X^{\underline{i}+p \underline{j} - p \delta_s}\,, \,\, {\rm since} \, 
i_s <p .
\end{array}
$$
This proves (b). 
\end{proof}

Now take $n=p^2$ (where $p$ is a prime) and consider the box partition
$\mathfrak b_p: p\geq p \geq \cdots \geq p$ of $n.$  
  Let $V^k_n:= k [\partial_{X_1}, \cdots,
\partial_{X_n}, \partial_{Y_1}, \cdots,
\partial_{Y_n}] (\triangle_{\mathfrak b_p}) \subset R_n^k$. 

Note that $\triangle_{\mathfrak b_p}$ has degree $<p$ in each variables $X_s$ and $Y_s$. The following 
lemma is clear.

\begin{lem} $\phi (\triangle_{\mathfrak b_p}) = \pm \triangle_n$ where $\triangle_n\in P_n^k$ is the
standard Vandermonde
determinant $\prod (X_s - X_t)$ (product being taken over $1 \leq s < 
t\leq n$) and $n:=p^2$.
\end{lem}

Let ${\mathbb D}_n^k$ be the algebra of all the divided differential operators 
in $n$ variables, i.e., the algebra over $k$ generated by $\{\partial_{X_s}^{(m)};
1\leq s \leq n, m\geq 0\}$ inside the algebra of all the linear operators of 
$P_n^k$. Since  $\partial_{X_s}^{(p^2)}$ kills
$\triangle_n$, by Lemma \ref{divided1}(b), we get the following. 

\begin{prop}  
$$\begin{array}{lll} \phi (V^k_n) &=& k [\partial_{X_1}, \cdots,
\partial_{X_n}, \partial_{X_1}^{(p)}, \cdots, \partial_{X_n}^{(p)}]
\triangle_n. \\
&=& {\mathbb D}_n^k (\triangle_n)\,,
\end{array}
$$ 
where $ {\mathbb D}_n^k (\triangle_n)$ denotes the space $\{D 
(\triangle_n): D\in {\mathbb D}_n^k\}$ inside  $P_n^k$.
 \end{prop}
\begin{cor}
\label{divided2} $\dim_k V^k_n \geq \dim_k {\mathbb D}_n^k (\triangle_n).$
\end{cor}

\begin{conjecture}
\label{conjecture3}
Let $n$ be any positive integer and let $p$ be a prime
such that $p^2\geq n$. Then, for any field $k$ of char. $p$, 
$$
\rm{dim}_k \,  {\mathbb D}_n^k (\triangle_n)= \rm{dim}_\mathbb C \,  
{\rm Diff}_n (\triangle_n),$$
where ${\rm Diff}_n (\Delta_n) :=\mathbb C[\partial_{X_1}, \cdots,
\partial_{X_n}] \cdot \Delta_n.$
\end{conjecture}

The latter of course is well known to have dimension $n!$.

\begin{remark} 
(a) It is easy to see that 
$$
\rm{dim}_k \,  {\mathbb D}_n^k (\triangle_n)\leq  \rm{dim}_\mathbb C \,  
{\rm Diff}_n (\triangle_n).$$

(b) It is possible that the above conjecture is true for any prime $p$
(not necessarily under the restriction  $p^2\geq n$), as some explicit 
calculations indicate, e.g., we have verified the conjecture for $
n=5, p=2$.

(c) The above conjecture is false, in general,  for the other polynomials which
transform under the action of the symmetric group $S_n$ via the sign character.
Take, e.g., $n=2$ and $f= (X_1-X_2)(X_1+X_2)$. Then, for $p=2$,
$\rm{dim}_k \,  {\mathbb D}_n^k (f)=2$, whereas $\rm{dim}_\mathbb C \,  
{\rm Diff}_n (f)=4$. 
 
\end{remark}

Assuming the above conjecture, we get the following theorem.

\begin{thm} Assume that $n=p^2$ and $k$ is a field of char. $p$. 
If the above Conjecture \ref{conjecture3} is true, then
$$ \rm{dim}_k \,V^k_n = n!.$$
In particular,
$$ \rm{dim}_\mathbb C \,{\rm Diff}_{2n}(\Delta_{\mathfrak b_p})= n!.$$

Thus the validity of Conjecture \ref{conjecture3} 
will imply the $n!$ result for the box 
partition $\mathfrak b_p$ of $n=p^2$.  
\end{thm}
\begin{proof} By combining Corollary \ref{divided2} and Conjecture 
\ref{conjecture3}, we get that
$ \rm{dim}_k \,V^k_n \geq n!.$ But, $ \rm{dim}_k \,V^k_n 
 \leq  \rm{dim}_\mathbb C \,{\rm Diff}_{2n}(\Delta_{\mathfrak b_p}) \leq n!.$ Thus,
 both the identities
 of the theorem follow. 
\end{proof}

\section{Gr Description of $A_\sigma$ (box plus one row case)}
\label{section8}
Let $n = pq+r$ with any integers $p> r \geq 0$ and $p>1, q\geq 1$, and 
consider the {\it box plus one row partition} $\sigma = 
\sigma(p,q,r) : \sigma_0 \geq \sigma_1 \geq \cdots \geq 
\sigma_q >0$ of 
$n$ with $\sigma_i = p$ for  $0 \leq i \leq q-1$ and 
$\sigma_q =r$. In this section we describe the 
algebra $A_{\sigma}$ as the associated graded ring of a 
specific filtration, depending only on $p$, of the 
cohomology ring $H^*(SL_n/B, \mathbb C)$. Here 
$B$ is the Borel subgroup of $SL_n$ consisting of 
the upper triangular matrices of determinant one. 
By a classical result we may identify  
$H^*(SL_n/B, \mathbb C)$ with the ring $P_n/J$, 
where $P_n :=\mathbb C[X_1, \cdots ,X_n]$ and 
$J \subset P_n$ denotes the ideal generated by the 
$S_n$-invariant elements of positive degree. 

\subsection{Preliminary results}
Recall the definitions of 
$S_{h,t,k}\in \mathbb C[Z_1, \cdots ,Z_t]$ 
and $n_k(\sigma)$ from the appendix.
In this section we will use the notation $n_k$ 
and $n_k^\vee$ to denote $n_k(\sigma)$ and 
$n_k(\sigma^\vee)$ respectively. Then, for the partition $\sigma = 
\sigma(p,q,r)$,   
\begin{equation}
\notag
n_k =
\begin{cases}
q(p-k) + (r-k) & \text{if $0 \leq k \leq r$}, \\
q(p-k) & \text{if $r \leq k < p$}, \\
0 & \text{if}\,\,k\geq p.
\end{cases}
\end{equation}
while 
\begin{equation}
\notag
n_k^\vee =
\begin{cases}
p(q-k) + r & \text{if $0 \leq k \leq q$}, \\
0 & \text{if}\,\, k>q.
\end{cases}
\end{equation}

Let $J_p$ denote the ideal in $P_n$ generated by $J$ and 
the elements $X_1^p, \dots, X_n^p$. Then  
\begin{lem}
\label{lem3}
Let $h$, $k$ and $1 \leq t \leq n$ be nonnegative integers.
If $h+t \geq n_k + 1$, then
$$S=S_{h,t,k}(X_{i_1},X_{i_2}, \dots, X_{i_t}) \in J_p,$$
whenever $1 \leq i_1 < i_2 < \dots < i_t \leq n$.
In particular, by Theorem \ref{procesi}, $\hat{J}(\sigma) \subset J_p$,
 where $\sigma = \sigma (p,q,r)$. 
\end{lem}
\begin{proof}
Observe first of all  that the statement is clear when
$k \geq p$. In the following we hence assume that $k <p$.
By symmetry it suffices to consider the case $i_j=j$,
$j=1,\dots,t$. Assume first
that $t \leq q$. Let $M=X_1^{s_1} X_2^{s_2} \cdots X_t^{s_t}$
be a monomial of degree $h$. If $s_j \leq p-k-1$ for
every $1\leq j \leq t$,  then
$$h= \operatorname{deg}(M) \leq t(p-k-1) \leq q(p-k)-t
\leq n_k - t\leq h-1,$$
which is a contradiction. Hence, there exists a $j$ such
that $s_j \geq p-k$. But then each term in $S$ is divisible by some $X_j^p$ 
and hence 
$S \in J_p$.

Assume now that $t \geq q+1$ and that the result, 
by induction, is true for smaller values of $t$. If
$k = 0$ then by Lemma 1.2 in \cite{ConProc} we find
$S \in J \subseteq J_p$. So assume that $k \geq 1$ 
and that the statement, by induction, is correct for
smaller values of $k$. As
$$(h+t) + t \geq h + t + q + 1  \geq n_k + q + 2
\geq n_{k-1} +1\,, $$
we know by induction that
$$S_{h+t, t, k-1}(X_1, \dots, X_t) \in J_p.$$
Rewriting $S_{h+t, t, k-1}(X_1, \dots, X_t)$ as
\begin{equation}
\notag
\begin{split}
& \sum_{1 \leq i_1 < i_2 < \cdots < i_l \leq t}
(X_1 \cdots X_t)^{k-1} (X_{i_1} X_{i_2}
\cdots X_{i_l}) S^l_{h+t-l}(X_{i_1},
\cdots, X_{i_l}) \\
= &
\sum_{1 \leq i_1 < i_2 < \cdots < i_l \leq t}
\frac{(X_1 \cdots X_t)^{k-1}}{(X_{i_1} \cdots
X_{i_l})^{k-1}} S_{h+t-l, l, k}(X_{i_1}, \cdots, 
X_{i_l}), \\
\end{split}
\end{equation}
and noticing, by induction, that all the terms with $l < t$
of the last expression belongs to $J_p$ , we conclude 
that $S_{h, t, k}(X_{1}, \cdots, X_{t}) \in J_p$.
\end{proof}

Let $J_q^\vee$ denote the ideal in $P_n$ generated 
by the elements of the form
$$(X_{i_1} \cdots X_{i_{r+1}})^q~,~1 \leq i_1 <
\cdots < i_{r+1} \leq n,$$
together with  the ideal $J_{q+1}$ (defined similar to $J_p$).

\begin{lem}
\label{lem4}
Let $h$, $k$ and $1 \leq t \leq n$ be nonnegative integers.
If $h+t \geq n_k^\vee + 1$ then
$$S=S_{h,t,k}(X_{i_1},X_{i_2}, \dots, X_{i_t}) \in J_q^\vee\,,$$
whenever $1 \leq i_1 < i_2 < \cdots < i_t \leq n$.
In particular, by Theorem \ref{procesi}, $\hat{J}(\sigma^\vee) \subset J^\vee
_q$,
where $\sigma = \sigma (p,q,r)$.
\end{lem}
\begin{proof}
Notice first of all that the statement is clear when
$k > q$. Hence, in  the following,  we  assume that $k \leq q$.
By symmetry it suffices to consider the case $i_j=j$, 
$j=1,\dots,t$. Assume first that $t \leq p$. 
 Assume, if possible, that $S
\notin J_q^\vee$. Then, by the definition of $J_q^\vee$, there exists a monomial
$M=X_1^{s_1} X_2^{s_2} \cdots X_t^{s_t}$ 
of degree $h+kt$ in the decomposition of $S$ such that 
at most $r$
of the indices $s_j$ are  $ \geq q$ and none of the indices are $\geq q+1$. Hence
$$h+tk= \operatorname{deg}(M) \leq rq + (t-r) (q-1)
\leq n_k^\vee +tk -t \leq h +tk-1,$$
which is a contradiction.

Thus we are left with the case $t > p$. Assume by
induction  that the
statement is correct for smaller values of $t$.
If $k = 0$ then by Lemma 1.2 in \cite{ConProc} we
find $S \in J \subseteq J_q^\vee$. So assume, by another
induction, that $k \geq 1$ and that the statement is
correct for smaller values of $k$. As
$$(h+t) + t \geq h + t + p + 1  \geq n_k^\vee + p + 2
\geq n_{k-1}^\vee +1 \,,$$
we know by induction that
$$S_{h+t, t, k-1}(X_1, \dots, X_t) \in J_q^\vee.$$
Using the identity (6) for $S_{h+t, t, k-1}(X_1, \dots, X_t)$ 
and noticing, by induction, that all the terms with $l<t$ 
belong to $J_q^\vee$, we conclude that
 $S_{h, t, k}(X_{1}, \cdots, X_t) \in J_q^\vee$.
\end{proof}

\begin{remark}
By the definition of $\hat{J}(\sigma)$ and $\hat{J}(\sigma^\vee)$ 
and Lemmas \ref{lem3}, \ref{lem4},  it follows easily that, in fact, 
  $J_p = \hat{J}
(\sigma)$ and $J_q^\vee = \hat{J}(\sigma^\vee)$.
\end{remark}

\vskip1ex

\subsection{The Gr description}

Let $x_i$ denote $X_i + J \in P_n/J$ and consider
the filtration of $P_n/J$:
$$F_0 \subseteq F_1 \subseteq F_2 \subseteq \cdots \,,$$
where $F_0 = \mathbb C$, $F_1 = \operatorname{span}_{\mathbb C}\{ 1,
x_1, x_2, \dots, x_n, x_1^p, x_2^p, \dots, x_n^p \}$
and $F_s$ is defined to be  $F_1^s$ (for any $s \geq 1$). Let $\operatorname{Gr}
(F)$
denote the associated graded algebra
$$\operatorname{Gr}(F) := \bigoplus_{s \geq 0}
F_s/F_{s-1},$$
where  $F_{-1} :=(0).$ 
Under the standard action of $S_n$ on $P_n$, each $F_s$ is $S_n$-stable and 
hence $\operatorname{Gr}(F)$ is canonically a $S_n$-module. 

Consider now the algebra homomorphism $\pi^\vee : P_n \rightarrow
\operatorname{Gr}(F)$ defined by $\pi^\vee(X_i)
= x_i + F_0 \in F_1/F_0$. By Lemma \ref{lem3},  $\hat{J}(\sigma)
\subseteq J_p \subseteq \operatorname{ker} (\pi^\vee)$
and thus, by Theorem \ref{procesi}, we get an induced map of graded algebras
$$\overline{\pi}^\vee : H^*(X_{\sigma^\vee}) \rightarrow
\operatorname{Gr}(F).$$
Consider also the algebra homomorphism $\pi
: P_n \rightarrow \operatorname{Gr}(F)$ defined
by $\pi(X_i) = x_i^p + F_0 \in F_1/F_0$.
As $p(q+1) \geq n$, we know that $X_i^{p(q+1)} \in J$
and hence $\pi(X_i^{q+1}) =0$.
Also, by  \cite{Haiman2}, Cor. 3.2.2,    
$$ (X_{i_1} \cdots X_{i_{r+1}})^{qp} \in J.$$
Thus, from Lemma \ref{lem4},   we conclude that $\hat{J} (\sigma^\vee)
\subseteq J_q^\vee \subseteq \operatorname{ker}
(\pi)$. This gives us an induced map
of graded algebras
$$\overline{\pi} : H^*(X_{\sigma})
\rightarrow \operatorname{Gr}(F).$$
Moreover, both $\overline{\pi}$ and $\overline{\pi}^\vee$ are 
$S_n$-equivariant. 

\begin{prop}
\label{gr2}
The map $\phi = \overline{\pi} \otimes \overline{\pi}^\vee
: H^*(X_{\sigma}) \otimes H^*(X_{\sigma^\vee})
\rightarrow \operatorname{Gr}(F)$ taking $a\otimes b \mapsto 
\overline{\pi} (a)\cdot \overline{\pi}^\vee (b)$ is an
$S_n$-equivariant  surjective graded
algebra homomorphism. Furthermore, the top graded component of $\operatorname{Gr}(F)$
is of degree $d_{\sigma}$ and is $1$-dimensional supporting the 
 sign representation of $S_n$.
\end{prop}
\begin{proof}
As $\operatorname{Im}(\phi)$, by definition, contains
$F_1/F_0$ and $F_0$ and as $\operatorname{Gr}(F)$
(by definition) is generated by degree one elements, we get
that $\phi$ is surjective. As the top graded component in
$H^*(X_{\sigma}) \otimes H^*(X_{\sigma^\vee})$
is of degree $d_{\sigma}$, the top  graded component in
$\operatorname{Gr}(F)$ is  of degree at most
 $d_{\sigma}$. As a representation of $S_n$,
the ring $\operatorname{Gr}(F)$ is of course isomorphic
to $P_n/J$ and thus   $\operatorname{Gr}(F)$ is isomorphic to the regular 
representation of $S_n$.
In particular,  $\operatorname{Gr}(F)$ contains a unique copy
of the sign representation of $S_n$. But
$H^*(X_{\sigma}) \otimes H^*(X_{\sigma^\vee})$
contains a unique copy of the sign representation
which occurs  in the top degree $d_{\sigma}$ (cf. the proof of Proposition 
\ref{our conjecture implies n!}).
Therefore, the copy of the sign representation
in
 $\operatorname{Gr}(F)$ must also be placed in
degree $d_{\sigma}$. So the proposition follows from the following. 
\end{proof}

\begin{prop}
The top graded component of $\operatorname{Gr}(F)$ is
$1$-dimensional.
\end{prop}
\begin{proof}
By Remark 3.4.(ii) in \cite{ConProc},  we know that the
top degrees of $H^*(X_{\sigma})$ and
$H^*(X_{\sigma^\vee})$ are respectively
$$ d = \frac{2qr + p (q-1) q}{2},$$
$$ d^\vee = \frac{qp(p-1) + r(r-1)}{2}.$$
 The image of the top degree in
$H^*(X_{\sigma}) \otimes H^*(X_{\sigma^\vee})$
under $\phi$, can therefore be represented by elements
in $P_n$ of degree 
$$ p d + d^\vee = \frac{n(n-1)}{2}.$$
But modulo $J$ there is only one element (up to 
a constant) in $P_n$ of this degree.
\end{proof}

\begin{cor} 
\label{gr description}
Let $n=pq+r$ be as in the beginning of this section. Then
there exists a $S_n$-equivariant isomorphism of  graded algebras:
$$\operatorname{Gr}(F) \simeq A_\sigma ,$$ 
where $\sigma = \sigma (p,q,r).$
\end{cor}
\begin{proof}
Define the ideal $I \subset \operatorname{Gr}(F)$ by
$$I=\{a\in \operatorname{Gr}(F): \text{the top degree component of}\, ab \,
\text{equals}\, 0\,\text{for all}\, b\in \operatorname{Gr}(F)\}.$$
Clearly $I$ is a $S_n$-stable graded ideal of $\operatorname{Gr}(F)$. Then
$\operatorname{Gr}(F)/I$ satisfies all the characterizing properties of 
Theorem \ref{characterize} (in view of Proposition \ref{gr2}). Thus 
$\operatorname{Gr}(F)/I \simeq T_\sigma$. But, by Theorem 
\ref{haiman algebra} and the $n!$ result,  we know that $T_\sigma$ 
has complex dimension $n! = {\rm dim}_\mathbb C 
\operatorname{Gr}(F)$. This forces $I=0$ and hence 
$\operatorname{Gr}(F) \simeq T_\sigma \simeq A_\sigma$.
\end{proof}

\begin{remark} (a) In the above proof we have used the $n!$ result for the partition $\sigma = \sigma (p,q,r)$. Conversely, the validity of the above corollary clearly implies the $n!$ result for $\sigma$.
\vskip1ex
(b) In \cite{Garsia-Haiman}, Section 6, a statement 
equivalent to Corollary \ref{gr description} for 2-row 
shapes (i.e., $q=1$) appears without proof. In the same 
paper it is noted that computational experiments suggest 
that Corollary \ref{gr description} should be true. Recently, M. 
Haiman  informed us that he had an unpublished proof of Corollary 
\ref{gr description}  (assuming the validity of the 
$n!$-conjecture for $\sigma(p,q,r)$) prior to our work. 
\end{remark}

\section{Appendix: Cohomology of Springer fibres}
\label{springer}

We recall some well known results on the cohomology of Springer fibres in the 
case of $G=SL_n(\mathbb C)$ (which we have used in the paper).

Let $\sigma : \sigma_0 \geq
\sigma_1 \geq \cdots \geq \sigma_m >0 $ be a partition of $n$ and let 
 $\sigma^\vee :
\sigma'_0 \geq \cdots \geq \sigma'_{m'}> 0 $ be the dual partition. For any
$m' < i \leq n$, set $\sigma'_i=0.$ For any $1 \leq k \leq n$, recall the 
definition of the integer $d_k(\sigma)$ from Section \ref{section2}, equation
 (\ref{dk}). 
 For any 
$k \geq 0$, let $n_k =n_k(\sigma) $ denote the integer
$$n_k := \sigma'_k + \sigma'_{k+1}  + \cdots .$$

 Let
$G=SL_n(\mathbb C)$ and let $X_\sigma$ be the Springer fibre (over complex numbers) corresponding to a fixed  
nilpotent matrix $M_\sigma$ (of size $n\times n$) with Jordan blocks of sizes 
$ \sigma_0,  \sigma_1, 
\cdots,  \sigma_m$. Let
 $H^*(X_{\sigma})$ denote
the cohomology ring of $X_{\sigma}$ with coefficients
in $\mathbb C$. 

The following description of $H^*(X_{\sigma})$ is due to
Tanisaki \cite{Tanisaki}, though in its  present form we have taken it from 
 \cite{Garsia-Procesi}, p. 84-85. The notation
$e_r$ denotes the $r$-th elementary symmetric polynomial. As earlier,  
$P_n :=  \mathbb C[X_1,\dots, X_n]$.  

\begin{thm}
\label{tanisaki}
As a graded algebra,  we have an $S_n$-equivariant
isomorphism
$$H^*(X_{\sigma}) \simeq P_n/J(\sigma),$$
where $J(\sigma)$ is the ideal generated by all
polynomials of the form
$$e_r(X_{s_1}, X_{s_2}, \dots, X_{s_k})$$
subject to the conditions
$1 \leq s_1 < \dots < s_k \leq n$ and all $k,r \geq 1$ satisfying
$k -d_k(\sigma) <r \leq k$.
\end{thm}

 \begin{defn}
Let $t \geq 1$ and $h,k \geq 0$ be integers.
Let $S_h^t \in \mathbb C[Z_1, \dots, Z_t]$ denote
the sum of all monomials  of degree $h$ in the variables $Z_1, \dots,
Z_t$. We then define
$$S_{h,t,k}= (Z_1 \dots Z_t)^k S^t_h
\in \mathbb C[Z_1, \dots,Z_t].$$
\end{defn}

 The following description of the
cohomology of $X_{\sigma}$ is taken from \cite{ConProc},  Theorem  2.2.

\begin{thm}
\label{procesi}
Let $\hat{J}(\sigma)$ denote the ideal in $P_n$ generated
by the polynomials $S_{h,t,k}(X_{i_1}, X_{i_2}, \dots, X_{i_t})$
for nonnegative integers $h$, $k$ and $1 \leq t \leq n$
subject to the condition $h + t = n_k +1$, and
$1 \leq i_1 < i_2 < \cdots < i_t \leq n$. Then
there exists an $S_n$-equivariant isomorphism
$$H^*(X_{\sigma^\vee}) \simeq P_n/\hat{J}(\sigma).$$
\end{thm}

\end{document}